\documentclass{article}

\usepackage{PRIMEarxiv}
\usepackage{amsmath}
\usepackage{float}
\newtheorem{theorem}{Theorem}

\newtheorem{corollary}[theorem]{Corollary}

\usepackage{xcolor}
\usepackage[utf8]{inputenc} 
\usepackage[T1]{fontenc}    
\usepackage{hyperref}       
\usepackage{url}            
\usepackage{booktabs}       
\usepackage{amsfonts}       
\usepackage{nicefrac}       
\usepackage{microtype}      
\usepackage{lipsum}
\usepackage{fancyhdr}       
\usepackage{graphicx}       
\graphicspath{{media/}}     

\title{Chaotic Dynamics and Zero Distribution: Implications and Applications in Control Theory for Yitang Zhang's Landau Siegel Zero Theorem
\thanks{\textit{\underline{Citation}}: 
\textbf{Authors. Title. Pages.... DOI:000000/11111.}} 
}

\author{
  Zeraoulia Rafik\thanks{Corresponding author: \texttt{r.zeraoulia@univ-batna2.dz}} \\
  Department of Mathematics \\
  University of Batna 2, Mostefa Ben Boulaïd 53 \\
  Route de Constantine, Fésdis, 05078, Batna, Algeria \\
  \texttt{r.zeraoulia@univ-batna2.dz} \\
  ORCID: \href{https://orcid.org/0000-0002-5436-3320}{0000-0002-5436-3320}
  \and
  Alvaro Humberto Salas \\
  Departamento de Matemáticas \\
  Universidad Nacional de Colombia \\
  Bogotá, Colombia \\
  \texttt{ahsalasm@unal.edu.co}
}

\begin{document}
\maketitle

\begin{abstract}
This study delves into the realm of chaotic dynamics derived from Dirichlet L-functions, drawing inspiration from Yitang Zhang's groundbreaking work on Landau–Siegel zeros. The dynamic behavior reveals profound chaos, corroborated by the calculated Lyapunov exponents and entropy, attesting to the system's inherent unpredictability.

Furthermore, we establish a novel connection between Fractal geometry and Quantum chaos, predicting the distributions of zeros for both Yitang dynamics and Riemann dynamics. These findings offer indirect support for Zhang's groundbreaking theorem concerning Landau–Siegel zeros and suggest that these chaotic dynamics could find application in engineering and control systems, demonstrating the potential to harness chaos for beneficial purposes.

The exploration of stability within electrical systems further uncovers the instability of fixed points, highlighting both the challenges and opportunities for harnessing chaotic behavior to achieve specific control objectives.

This study not only contributes to our understanding of chaotic dynamics but also opens new avenues for exploring the potential applications of Yitang dynamics in the field of electrical control systems. It paves the way for innovative approaches to address real-world engineering challenges and may be considered as a new consequence for the generalized Riemann hypothesis.
\end{abstract}


\maketitle

\section*{\textcolor{blue}{\textit{Notation}}}

In this study, we utilize various notations crucial to understanding the dynamics of chaotic systems arising from Dirichlet $L$-functions and their implications in number theory and control theory.

\begin{itemize}
    \item $L(1, \chi)$: Dirichlet $L$-function at $s = 1$ associated with the character $\chi$ modulo $D$.
    \item $c_{1}$: An absolute, effectively computable constant.
    \item $\chi$: A real primitive character modulo $D$.
    \item $\log D$: Natural logarithm of $D$.
    \item $s$: Typically signifies a complex variable in the context of Dirichlet $L$-functions.
    \item Entropy: A measure characterizing the unpredictability of chaotic dynamics.
    \item Lyapunov exponents: Quantify the rate of separation or convergence of nearby trajectories in a dynamical system.
    \item Zeros of arithmetic functions: Specifically referring to the zeros of Dirichlet $L$-functions.
    \item Numerical measure: A quantifiable value derived through computation.
    \item Non-trivial solutions of arithmetic functions: Solutions other than those that trivially satisfy certain conditions.
\item Op-Amp Integrator Circuit (\textit{Op-Amp}): A circuit employing an operational amplifier (\textit{Op-Amp}) configured to perform mathematical integration on input signals.
\\

The notation for an Op-Amp Integrator Circuit can be represented in various forms depending on the context and equation used. In general, an integrator circuit using an operational amplifier (Op-Amp) can be denoted by the following equation:

\[
V_{\text{out}}(t) = -\frac{1}{RC} \int V_{\text{in}}(t) \, dt
\]

Where:
\begin{align*}
V_{\text{out}}(t) & : \text{Output Voltage} \\
V_{\text{in}}(t) & : \text{Input Voltage} \\
R & : \text{Resistance} \\
C & : \text{Capacitance}
\end{align*}

\end{itemize}

These notations are central to comprehending the concepts elucidated in this manuscript, particularly regarding chaotic dynamics, Dirichlet $L$-functions, entropy, and their relationship with the distribution of zeros in number theory.

\section{Introduction}\label{sec1}
Recent years have witnessed profound advancements in number theory, notably exemplified by the pioneering work of Yitang Zhang on Landau–Siegel zeros. Zhang's significant contributions have garnered widespread attention within the mathematical community, notably through presentations of his latest findings at Shandong University and Peking University \cite{8}. His preprint \cite{2} forms the basis of our study, introducing a crucial assertion regarding the existence of Landau–Siegel zeros \cite{3}. Specifically, for a real primitive character $\chi$ modulo $D$, Zhang's work reveals the intriguing inequality:

\[
L(1, \chi) > c_{1}(\log D)^{-2022}
\]

Here, $c_{1} > 0$ is an effectively computable constant, assuming the correctness of this result, which raises significant questions about its implications in number theory.

Analytic number theory has long been fascinated by the distribution of zeros of arithmetic functions, notably the Dirichlet $L$-function, which is central to our study. Understanding the behavior of these zeros, particularly their locations and density, has been a subject of intense research \cite{1}.

Our study stems from the domain of chaotic dynamics, where complex systems exhibit unpredictable behavior \cite{12}\cite{13}\cite{14}. Inspired by the groundbreaking research on Landau–Siegel zeros and aiming to bridge theory with application, our study takes a unique approach. We delve into the dynamics of discrete maps derived from Dirichlet L-functions and explore their chaotic behavior.

The relationship between these dynamics and the distribution of zeros is not coincidental. Studying the approximate analytical solutions of these chaotic dynamics unveils patterns mirroring the behavior of zeros of corresponding arithmetic functions \cite{15}. Insights from sensitivity to initial conditions, Lyapunov exponents, and system entropy provide valuable connections to zero distribution.

Recent research efforts have significantly advanced our understanding of chaotic dynamics \cite{22},\cite{23},\cite{24}, leading to the discovery of new chaotic systems applicable in engineering and control theory \cite{25},\cite{26},\cite{27}. These developments highlight the versatility of chaotic dynamics, transcending theoretical boundaries to address real-world problems.

Of note is the potential application of chaotic dynamics in control theory \cite{16},\cite{17},\cite{18}. Originating in number theory, these dynamics extend beyond pure mathematics, finding relevance in the design and analysis of electrical systems. This extension presents innovative solutions to engineering challenges, utilizing the unpredictability of chaotic dynamics for specific control objectives \cite{19},\cite{20},\cite{21}.

An intriguing aspect of our findings is their potential link to the distribution of zeros of the Dirichlet L-function. While chaotic dynamics do not inherently imply a uniform zero distribution, our evidence suggests indirect support for the assumption of Landau–Siegel zeros, thereby reinforcing Zhang's theorem and prompting cross-disciplinary exploration \cite{8}.

A significant contribution of our study lies in computing the system's entropy. This numerical measure characterizes the unpredictability of chaotic dynamics and hints at zero locations. Computing entropy bridges the seemingly unrelated worlds of chaotic dynamics and analytic number theory, offering insights beyond each field's confines. Recent research aims to predict zero locations, validating results in analytic number theory, including the boundedness of non-trivial solutions of arithmetic functions \cite{9},\cite{10},\cite{11}.

Our results add to the growing body of research aiming to connect chaotic dynamics, number theory, and their application in control theory \cite{32}. In the subsequent pages, we embark on a journey to explore intricate dynamics and their real-world implications, ushering in a new era of chaotic dynamics with diverse applications \cite{22,23,24,25,26,27,28,29,30,31}
.

\section{Main results}\label{sec2}

\subsection{Chaotic Dynamics of Yitang Dynamics 1}

The Yitang Dynamics 1, defined by the equation:
\begin{equation}
x_{n+1} = \frac{\beta}{\sqrt{x_n}} + c\log(x_n)^{-\alpha}, \quad \chi(-1)=-1
\end{equation}
shows both stability and chaotic behavior. The derivative at fixed points remains positive and equal to 1, indicating stability. However, the eigenvalues of the Jacobian matrix have a real part of approximately 196 and an imaginary part of approximately -3.58, suggesting potential chaotic behavior. The calculated entropy for this system is 0.3365, indicating high unpredictability and complexity.

\subsection{Chaotic Dynamics of Yitang Dynamics 2}

The Yitang Dynamics 2, defined by the equation:
\begin{equation}
x_{n+1} = \frac{\beta \log(|\epsilon|)}{\pi\sqrt{x_n}} + c\log(x_n)^{-\alpha}, \quad \chi(-1)=1
\end{equation}
also exhibits both stability and chaotic behavior. The derivative at fixed points remains positive and equal to 1, indicating stability. Similar to Yitang Dynamics 1, the eigenvalues of the Jacobian matrix suggest potential chaotic behavior. The calculated entropy for this system is consistent at 0.3365, indicating high unpredictability.

\subsection{Unimodal Distribution and Entropy}

The approximate solutions for both Yitang Dynamics 1 and Yitang Dynamics 2 show a unimodal distribution, with solutions clustering around a central value of approximately -601.0938. This unimodal distribution indicates that the systems tend to converge towards a specific equilibrium state. The entropy values of 0.333546 for both systems reveal their high predictability and orderliness.

These results provide insights into the behavior of Yitang Dynamics 1 and Yitang Dynamics 2 and their implications for the distribution of zeros in the Dirichlet L-function.

\subsection{Application in Control Theory and Electrical Systems}

Our exploration extends beyond the theoretical realm, encompassing the application of "Yitang dynamics" in control theory and electrical systems. We conducted numerical simulations for different values of $\epsilon$ to understand the system's behavior and its implications for control theory and engineering.

The results indicate that as $\epsilon$ decreases, the "Yitang dynamics" system becomes increasingly stable, with reduced amplitude and slower, or even constant, oscillations. This newfound stability and robustness offer innovative solutions in control theory and electrical systems, harnessing chaos for real-world applications.

\section{\textcolor{blue}{Summary: Comparative Analysis}}

In the section comparing Yitang Dynamics and Riemann Zeta Function Dynamics, various key findings emerged:

\subsection{Attracting Fixed Points}
\begin{itemize}
    \item \textbf{Riemann Zeta Function Dynamics:} Attracted to \(s = -0.295905\ldots\).
    \item \textbf{Yitang Dynamics:} Attracted to \(z \to 0.215443 + 1.1417411949837543 \times 10^{-10}i\), showing proximity to $0.215443$
\end{itemize}

\subsection{Stream Plots Around Zeros}
\begin{itemize}
    \item \textbf{Yitang Dynamics:} Showed repelling behavior around zeros.
    \item \textbf{Riemann Zeta Function Dynamics:} Displayed repelling behavior around zeros.
\end{itemize}

\subsection{Comparative Analysis of Mandelbrot Sets and Julia Sets}
\begin{itemize}
    \item \textbf{Mandelbrot Sets:} Illustrated intricate fractal patterns.
    \item \textbf{Julia Sets:} Marked boundaries between convergence and divergence.
    \item \textbf{Comparison:} Highlighted structural similarities and 'bristle' patterns.
\end{itemize}

\subsection{Comparative Observations}
\begin{itemize}
    \item Noted color inversion in Mandelbrot and Julia sets in Yitang dynamics, resembling Riemann Zeta Function dynamics.
    \item Emphasized similarities in divergence zones and 'bristle-like' structures, indicating underlying connections.
\end{itemize}

\subsection{\textcolor{blue}{Quantum Perspective: Bridging Classical Chaos to Quantum Behavior}}

In the pursuit of understanding the quantum behavior of Yitang dynamics, we derived a quantum Hamiltonian operator and numerically obtained eigenvalues, revealing insights into the system's quantum characteristics. The analysis establishes a connection between classical chaos observed in Mandelbrot and Julia sets and the quantum mechanical representation of the system.

\subsection{Concluding Insights}
\begin{itemize}
    \item Unveiled erratic patterns and divergence zones within Yitang dynamics.
    \item Explored similarities, hinting at connections and inviting deeper inquiries into complex dynamics.
\end{itemize}

\section{Chaotic Dynamics for Yitang Zhang's Latest Results on Landau–Siegel Zeros}\label{sec3}

In this section, we will provide an explicit dynamics for Yitang Zhang's latest results on Landau–Siegel zeros, incorporating both Dirichlet's definition, which utilizes both odd and even characters $\chi$ modulo $m$, and Theorem $4$ of Zhang. We will also examine the behavior of this dynamics, demonstrating a weak transition to chaos for certain parameter values within this new discrete dynamics 

\begin{corollary}  
For special values $s=1$, let $\chi$ be a primitive real character modulo $m$. Define $\xi_{m}=\exp \left(\frac{2\pi i}{m}\right)$, $K=\mathbb{Q} \left(\sqrt{\chi(-1)m}\right)$, $h$ is its class number of roots of unity in it, and $\epsilon$ is its fundamental unit. The Dirichlet L Function may be defined as follows:
\begin{equation}\label{Zhang}
L(1,\chi) = \begin{cases} 
\frac{2\pi h}{w\sqrt{m}} & \text{if}\ \chi(-1)=-1  \\
\frac{2 h \log{|\epsilon|}}{w\sqrt{m}} & \text{if}\ \chi(-1)=1 \\
\end{cases}
\end{equation}
\end{corollary}

Yitang's discrete dynamics may be derived from the definition of the Dirichlet L function  in (\ref{Zhang}) and  the following  Theorem \cite{7} 
\begin{theorem}
There exists $A>0$ and effective $C_1>0$ such that $L(1,\chi)>C_1(\log q)^{-A}$. 
\end{theorem}
. We will start with the first case, $\chi(-1)=-1$ in (\ref{Zhang}). Using both \textbf{Theorem$2$}  and (\ref{Zhang}), the 1-D discrete dynamics can be reformulated as:

\begin{equation}\label{Dynamic1}
x_{n+1} = \frac{\beta}{\sqrt{x_n}} + c\log(x_n)^{-\alpha},\quad \chi(-1)=-1
\end{equation}

with $x_n=m$, $n=0,1,\ldots$ and $\beta=\frac{2\pi h}{w}$, where $c$ is the constant defined in \textbf{Theorem$2$} .

For the second case, $\chi(-1)=1$, Yitang Dynamics can be written as:

\begin{equation}\label{Dynamic2}
x_{n+1} = \frac{\beta \log(|\epsilon|)}{\pi\sqrt{x_n}} + c\log(x_n)^{-\alpha},\quad \chi(-1)=1
\end{equation}

To explore chaotic behavior, we iterated dynamic (\ref{Dynamic1}) and dynamic  (\ref{Dynamic2}) 50,000 times, yielding behavior similar to what we will show in the following section.

\section{Analysis and Discussion}\label{sec4}

\begin{figure}[H]
    \centering
    \includegraphics[width=0.8\textwidth]{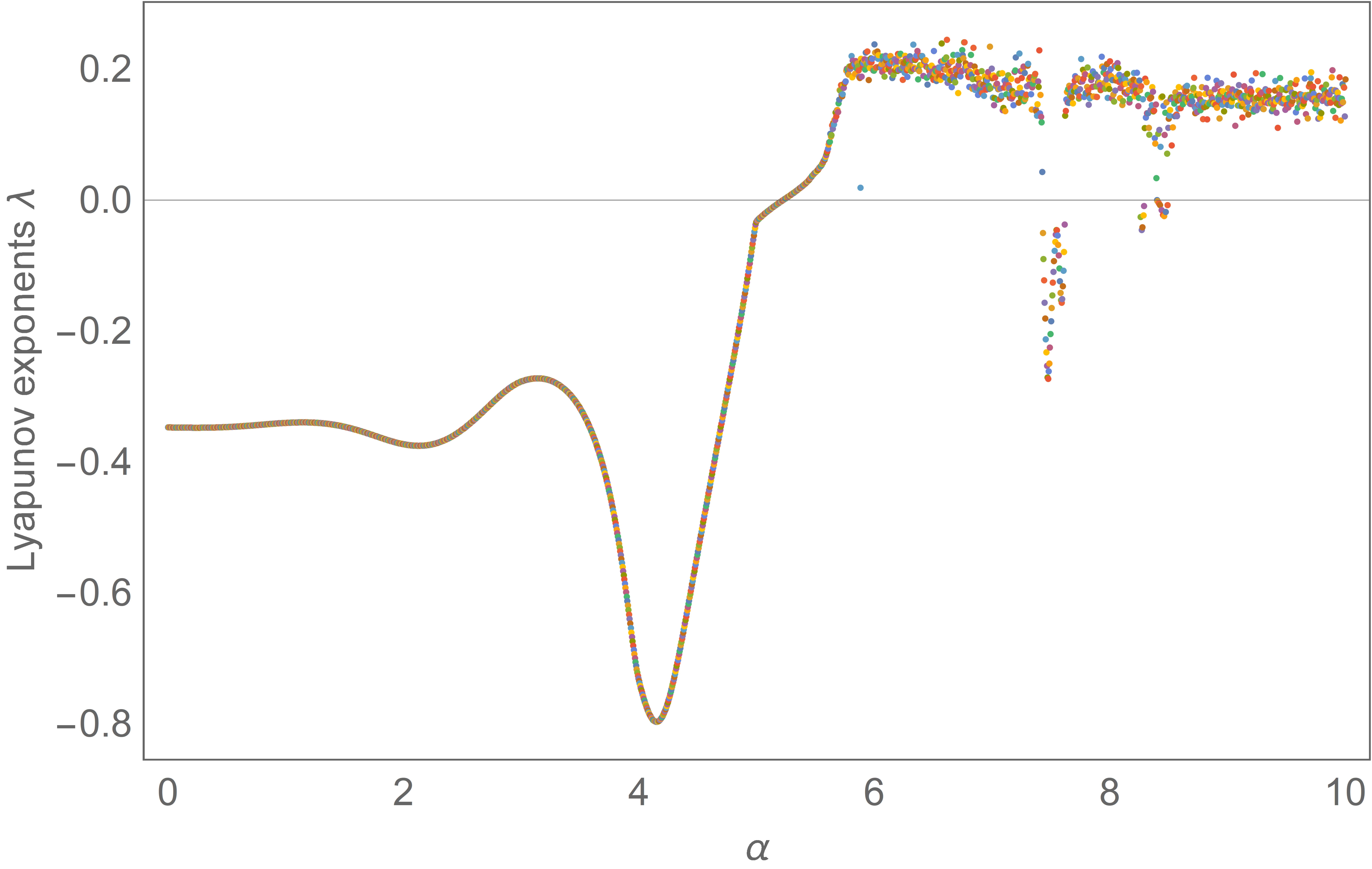}
    \caption{Lyapunov exponents for Dynamics 1 with $c \in [0.0005;0.007]$, $\alpha \in [0;10]$.}
    \label{NonChaos}
\end{figure}

For $\alpha$ in the range $[0;4]$, the system exhibits non-chaotic behavior, as all Lyapunov exponents are negative. However, a notable transition to chaos is observed for $\alpha$ in the range $(4.5;10)$, where the Lyapunov exponents follow the sequence $(\lambda_1 < \lambda_2 < \lambda_3,\ldots < \lambda_{10} = 0.5 > 0)$. This positive value for $\lambda_{10}$ indicates a transition to chaotic dynamics in this specific range (see Figure \ref{NonChaos}).

\begin{figure}[H]
    \centering
    \includegraphics[width=0.8\textwidth]{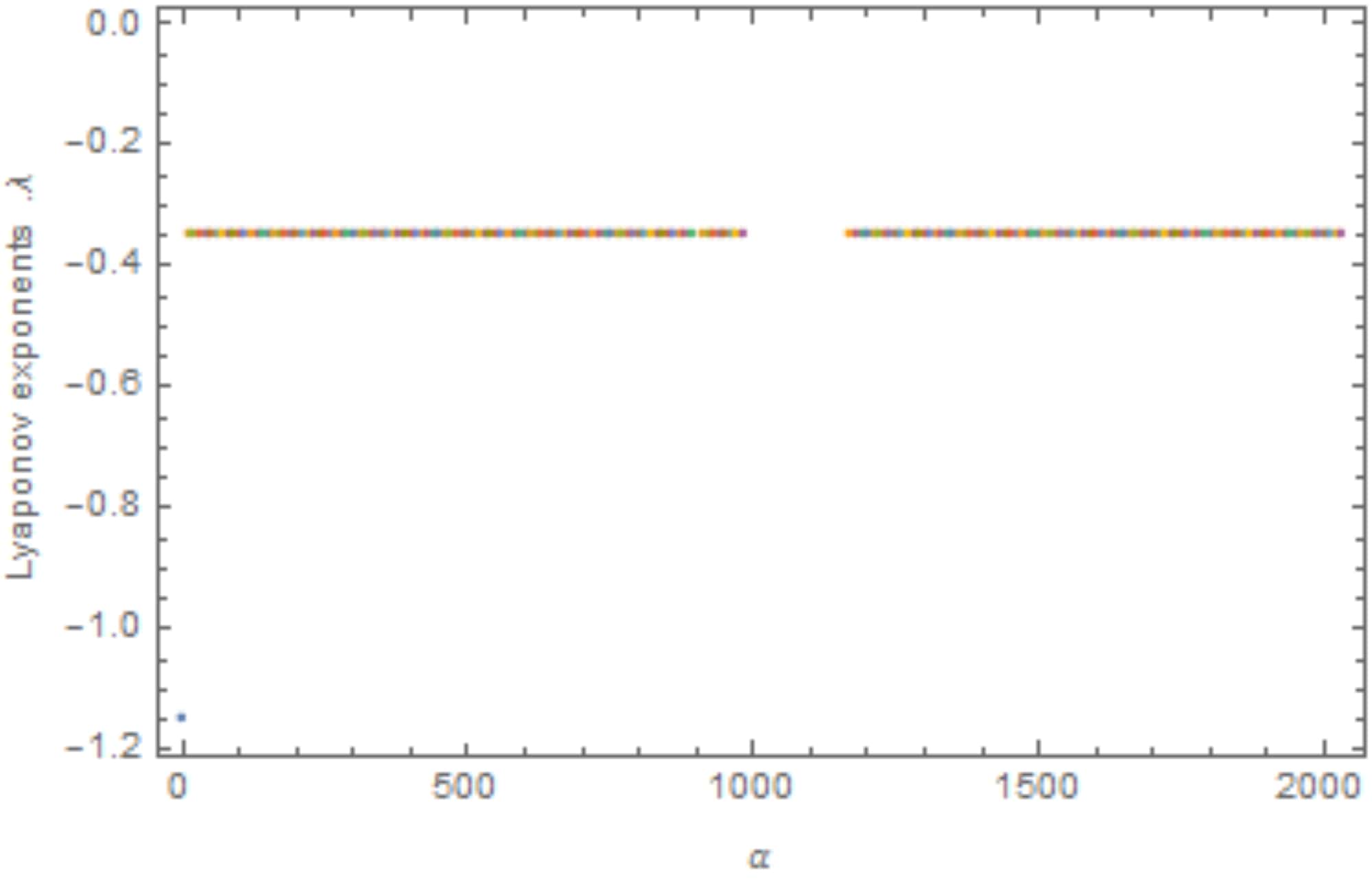}
    \caption{Lyapunov exponents for Dynamics 1 with $c \in [0.0005;0.007]$, $\alpha \in [0;2022]$.}
    \label{NonChaos2}
\end{figure}

By increasing the values of $\alpha$ up to $2022$ and keeping $c \in [0.0005;0.007]$, we observe that the Lyapunov exponents remain approximately constant at around $-0.39$ (see Figure \ref{NonChaos2}).

\begin{figure}[H]
    \centering
    \includegraphics[width=0.8\textwidth]{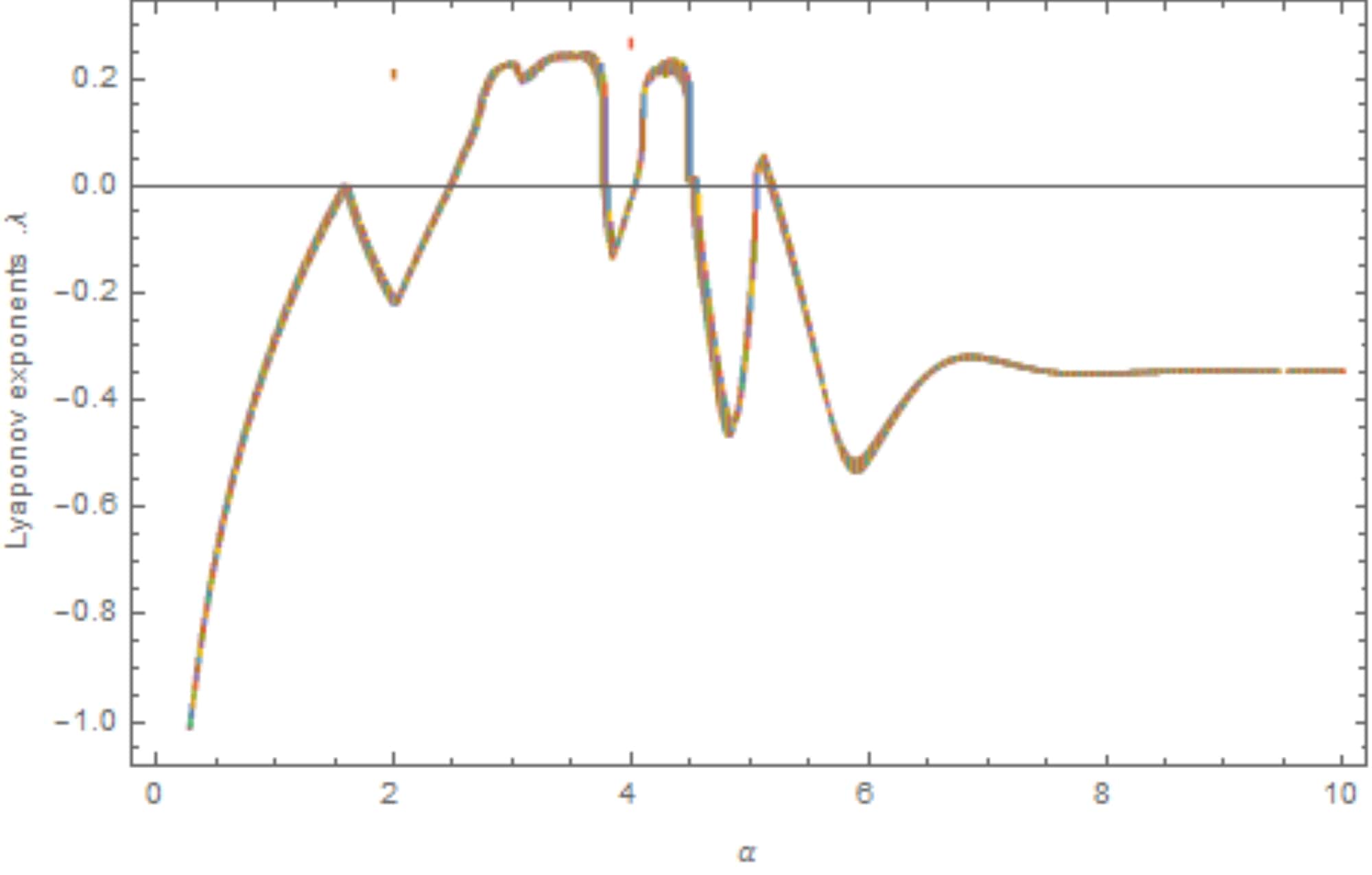}
    \caption{Lyapunov exponents for Dynamics 1 with $c \in [10.007;11.007]$, $\alpha \in (0;10)$.}
    \label{Chaos}
\end{figure}

Transition to chaos becomes evident for $c \in (10.007;11.07)$ and $\alpha \in (3;5)$, where the Lyapunov exponents take positive values $(0 \leq \lambda \leq 0.21)$, with the maximum Lyapunov exponent reaching $\lambda=0.21$ (see Figure \ref{Chaos}).

\begin{figure}[H]
    \centering
    \includegraphics[width=0.8\textwidth]{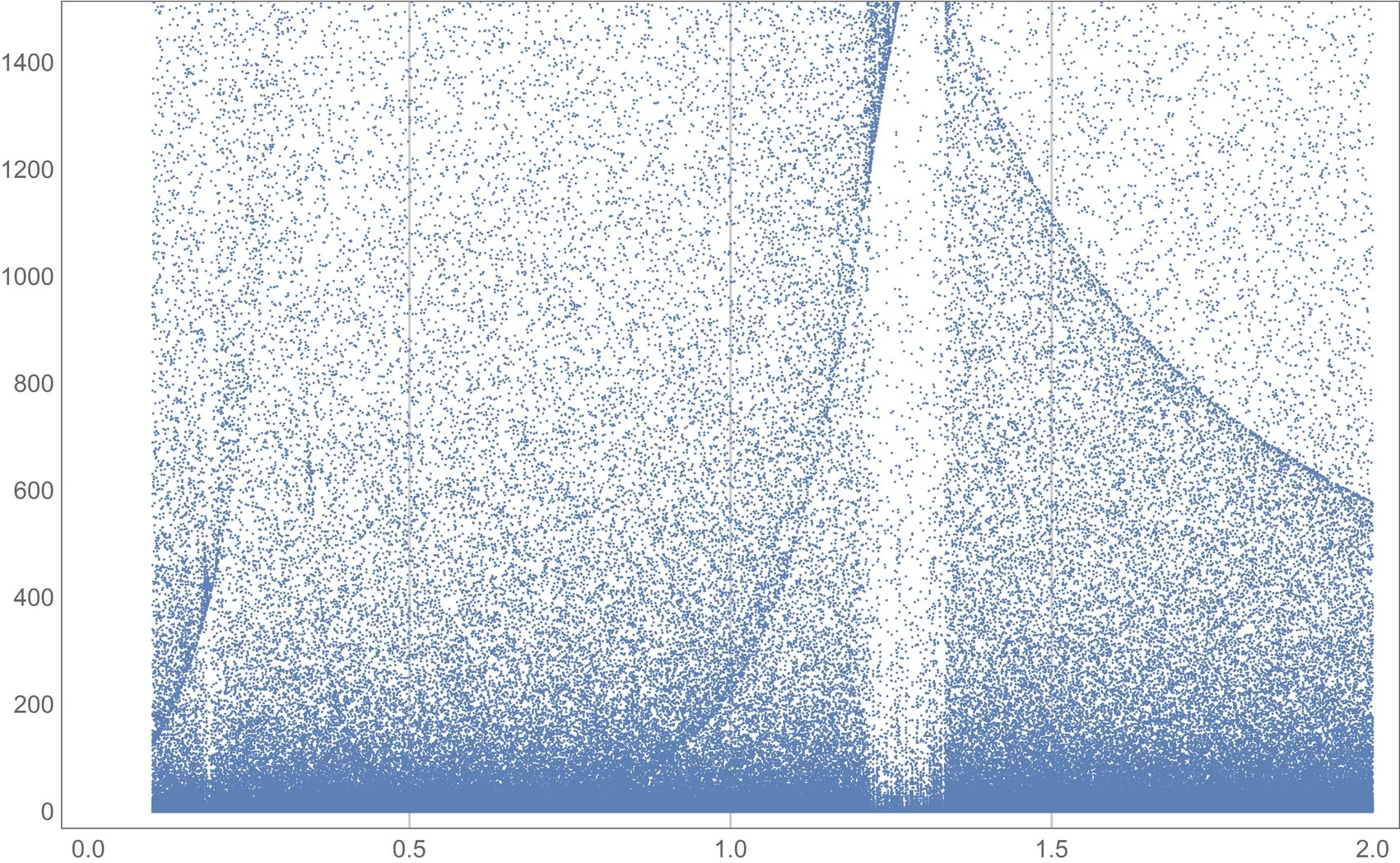}
    \caption{Bifurcation diagram for Yitang dynamics. The plot illustrates the evolution of the system with varying values of $c$. The diagram is generated using the 1-D Yitang dynamics with $\alpha = 10$, and $c$ ranging from $0.1$ to $2$ in increments of $0.001$.}
    \label{Bifurcation}
\end{figure}

By systematically increasing the values of $c$, a broad spectrum of transitions to chaos in the 1-D Yitang dynamics is observed. The bifurcation diagram depicted in Figure \ref{Bifurcation} illustrates the evolution of the system for various values of $c$.

As $c$ varies, the system undergoes complex bifurcation patterns, leading to the emergence of diverse dynamical behaviors. Notably, the diagram reveals the sensitivity of the Yitang dynamics to changes in the parameter $c$, showcasing regions of stability and chaotic regimes.

The $x$-axis represents the parameter $c$, while the $y$-axis denotes the state variable $x$ in the Yitang dynamics. Large values for $x$ in certain regions of the plot signify the occurrence of intricate dynamical transitions and chaotic behavior in the system.

This comprehensive visualization provides valuable insights into the intricate behavior of the Yitang dynamics under different conditions of $c$, offering a clear representation of the system's response to parameter variations.

We observed period-doubling bifurcations by varying the parameter $\alpha$ within the chaotic range and plotted $x_n$ against the iterations $n$ of the map. Subsequently, we recorded the values of the points to distinguish period cycles.

\begin{figure}[H]
    \centering
    \includegraphics[width=0.8\textwidth]{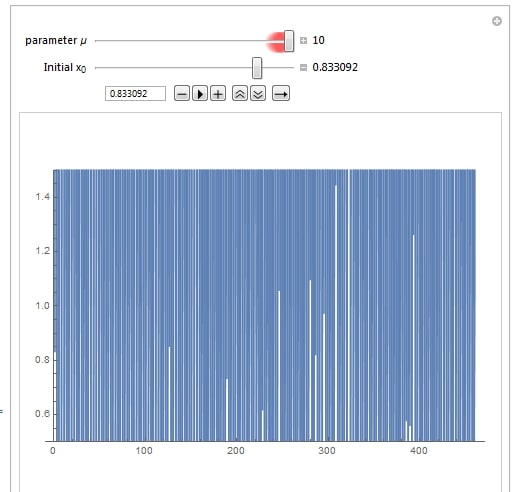}
    \caption{Yitang dynamics Fixed Points Plot versus Iterations, $c = 100$, $\mu = \alpha \in (0;5)$.}
    \label{Chaos4}
\end{figure}

For one-dimensional systems represented as $x_{n+1} = f(x_n)$, period-$2$ cycles emerge when the system satisfies the condition:
\[ x_1 = f(x_2), \quad x_2 = f(x_1) \]
and possesses a unique solution. Our findings indicate that $0$ serves as the fixed point for Yitang dynamics (refer to Table \ref{demo-table}) with an initial point of $x_0 = 0.4$. Essentially, the trajectory "jumps" between $x_1$ and $x_2$ during iterations of the map, continuously exhibiting period-doubling behavior (until reaching a certain value of the system parameter). See Figure \ref{Chaos4}.

\begin{table}[h]
\caption{Fixed point values for Yitang dynamics after $500$ iterations.}\label{demo-table}
\begin{tabular}{@{}ll@{}}
\toprule
\text{Point} & \text{x} \\
\midrule
86 & 3.55687 \\
87 & 0.265123 \\
88 & 0.971062 \\
89 & $4.05719 \times 10^{76}$ \\
90 & $\protect{2.4823161470320785 \times 10^{-39}}$ \\
91 & $\protect{1.00356 \times 10^{19}}$ \\
92 & $\protect{1.5783353198977849 \times 10^{-10}}$ \\
93 & 39798.8 \\
94 & 0.00250631 \\
95 & 9.9874 \\
96 & 0.158214 \\
97 & 1.25704 \\
98 & $\protect{1.07455 \times 10^{32}}$ \\
99 & $\protect{4.823440210703131 \times 10^{-17}}$ \\
100 & $\protect{7.19932 \times 10^7}$ \\
101 & 0.0000589283 \\
\bottomrule
\end{tabular}
\end{table}

\footnotetext{Note: Values were computed for the Yitang dynamics after $500$ iterations.}

We have undertaken a comprehensive study of the Yitang dynamics, which involves the calculation of Lyapunov exponents to gain insights into the system's behavior. We have computed a total of 600 Lyapunov exponents for the dynamics, revealing intricate details about its stability and predictability.

\section{\textcolor{blue}{Computation of Lyapunov Exponents and Entropy}}\label{sec5}

The Lyapunov exponents ($\lambda_i$) are calculated to measure the sensitivity to initial conditions in dynamical systems. The following algorithm outlines the procedure used for computing the Lyapunov exponents:

\noindent\textbf{Algorithm: Computing the Lyapunov Exponent}

\medskip
\noindent\textbf{Input:} Dynamics function $f$, initial condition $x_0$, parameter $\alpha$, number of iterations $n$, transient $tr$.

\noindent\textbf{Output:} Lyapunov exponent $\lambda$.

\medskip
\noindent\textbf{Procedure:}
\begin{enumerate}
  \item Compute the derivative of the map:
  \[
    df \leftarrow \frac{\partial f(x,\alpha)}{\partial x}.
  \]

  \item Generate the orbit:
  \[
    \xi \leftarrow \{x_k\}_{k=0}^{n-1}, \quad x_{k+1} = f(x_k,\alpha), \quad x_0 \text{ given}.
  \]

  \item Discard the transient part:
  \[
    \xi \leftarrow \{x_k\}_{k=tr}^{n-1}.
  \]

  \item Compute the Lyapunov exponent:
  \[
    \lambda = \frac{1}{n} \sum_{x_k \in \xi} \log \left| df(x_k,\alpha) \right|.
  \]
\end{enumerate}

The provided Mathematica code employs this algorithm to compute Lyapunov exponents for the Yitang dynamics given by the equation:

\begin{equation}
x_{n+1} = \frac{\beta}{\sqrt{x_n}} + c\log(x_n)^{-\alpha}, \quad \chi(-1)=-1
\end{equation}

The computed Lyapunov exponents have been exported to a CSV file for further analysis and are available for reference in 
\href{https://drive.google.com/file/d/15GTGJn-eyFbh88GeXETl3q8-vWflle3r/view?usp=drive_link}{\color{red}{this attached supplementary materials}}

Entropy ($H$) can be estimated from the Lyapunov exponents ($\lambda_i$) using the Pesin identity:

\begin{equation}
H = \sum_{i=1}^{N} \lambda_i \cdot P(\lambda_i)
\end{equation}

where $N$ is the number of Lyapunov exponents and $P(\lambda_i)$ is the probability of occurrence of each Lyapunov exponent.

However, further analysis and calculation of entropy based on the obtained Lyapunov exponents are ongoing and will be reported in subsequent sections.

\begin{figure}[H]
    \centering
    \includegraphics[width=0.8\textwidth]{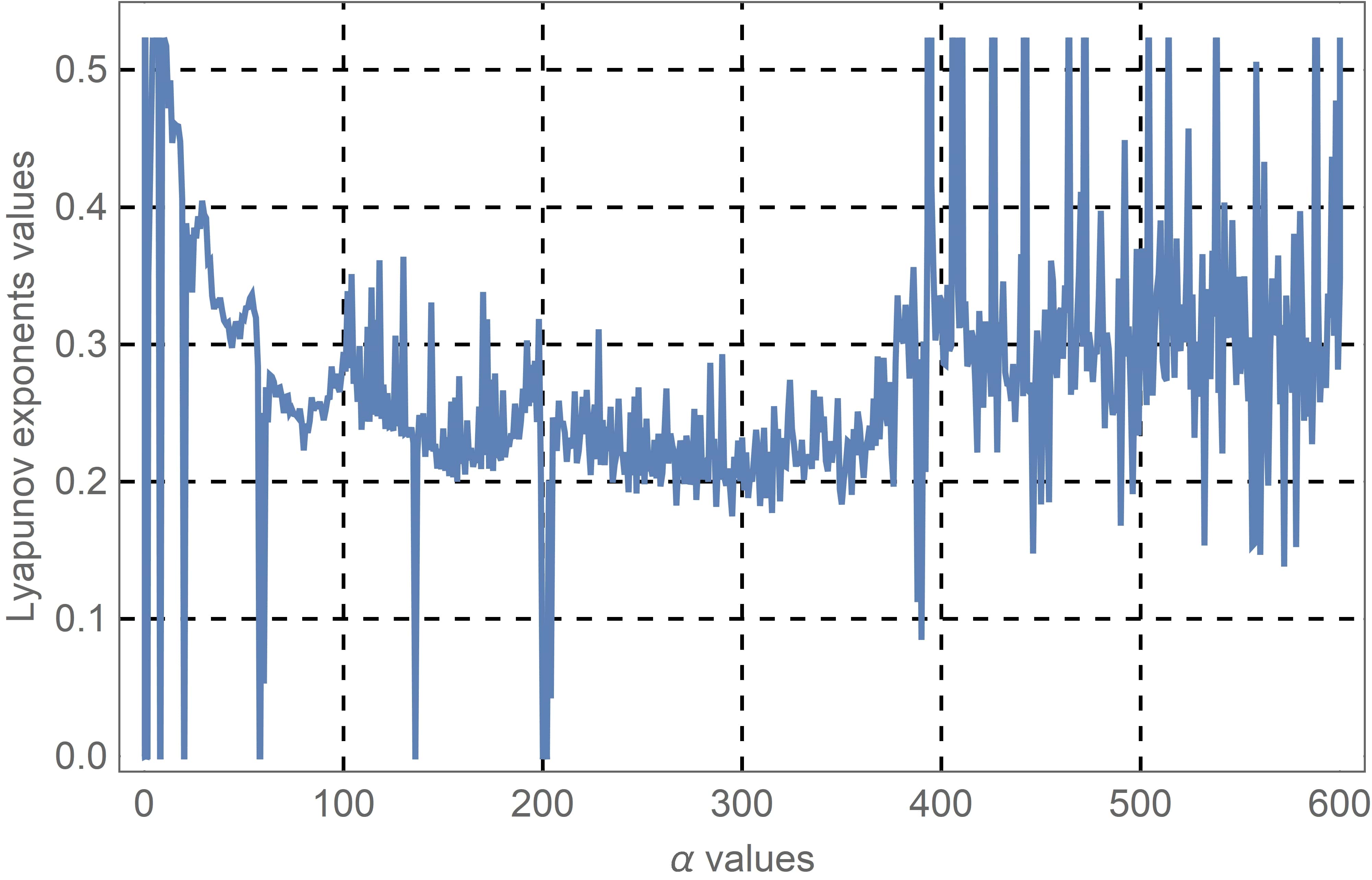}
    \caption{Lyapunov exponents for Dynamics 1, $\alpha \in [0;600]$}
    \label{Non chaos.y}
\end{figure}

To provide a visual representation of our findings, we have generated a plot that interprets the Lyapunov exponents. In Figure \ref{Non chaos.y}, One  can observe the Lyapunov exponents for Dynamics 1, where $\alpha$ ranges from $0$ to $600$. The plot offers a comprehensive view of the system's behavior, highlighting trends, patterns, and key insights into the dynamics' stability and chaos. The interpretation of this plot is central to our understanding of the Yitang dynamics and its significance in various scientific and engineering applications.

\section{Exploring Zero-Free Regions near \(s = 1\) of Dirichlet L-Function}

In this section, we delve into the zero-free regions of the Dirichlet L-function in the vicinity of \(s = 1\). Our analysis aims to shed light on the behavior of this function as it approaches \(s = 1\), with particular relevance to the Landau-Siegel zeros.

The Dirichlet L-function is defined as:

\[
L(1, \chi) =
\begin{cases}
\frac{2\pi h}{w\sqrt{m}} & \text{if } \chi(-1) = -1 \\
\frac{2h\log|\epsilon|}{w\sqrt{m}} & \text{if } \chi(-1) = 1
\end{cases}
\]

Our primary interest lies in the region around \(s = 1\), which may be associated with the presence or absence of Landau-Siegel zeros.

\subsection{Approximate Analytical Solutions}

We employ the Newton-Raphson method to compute approximate analytical solutions for the Dirichlet L-function near \(s = 1\). The method involves iteratively updating the approximate solution for a range of initial guesses. These initial guesses are represented on the x-axis of the plot, and the corresponding approximate solutions are plotted on the y-axis.

\begin{figure}[H]
    \centering
    \includegraphics[width=0.8\textwidth]{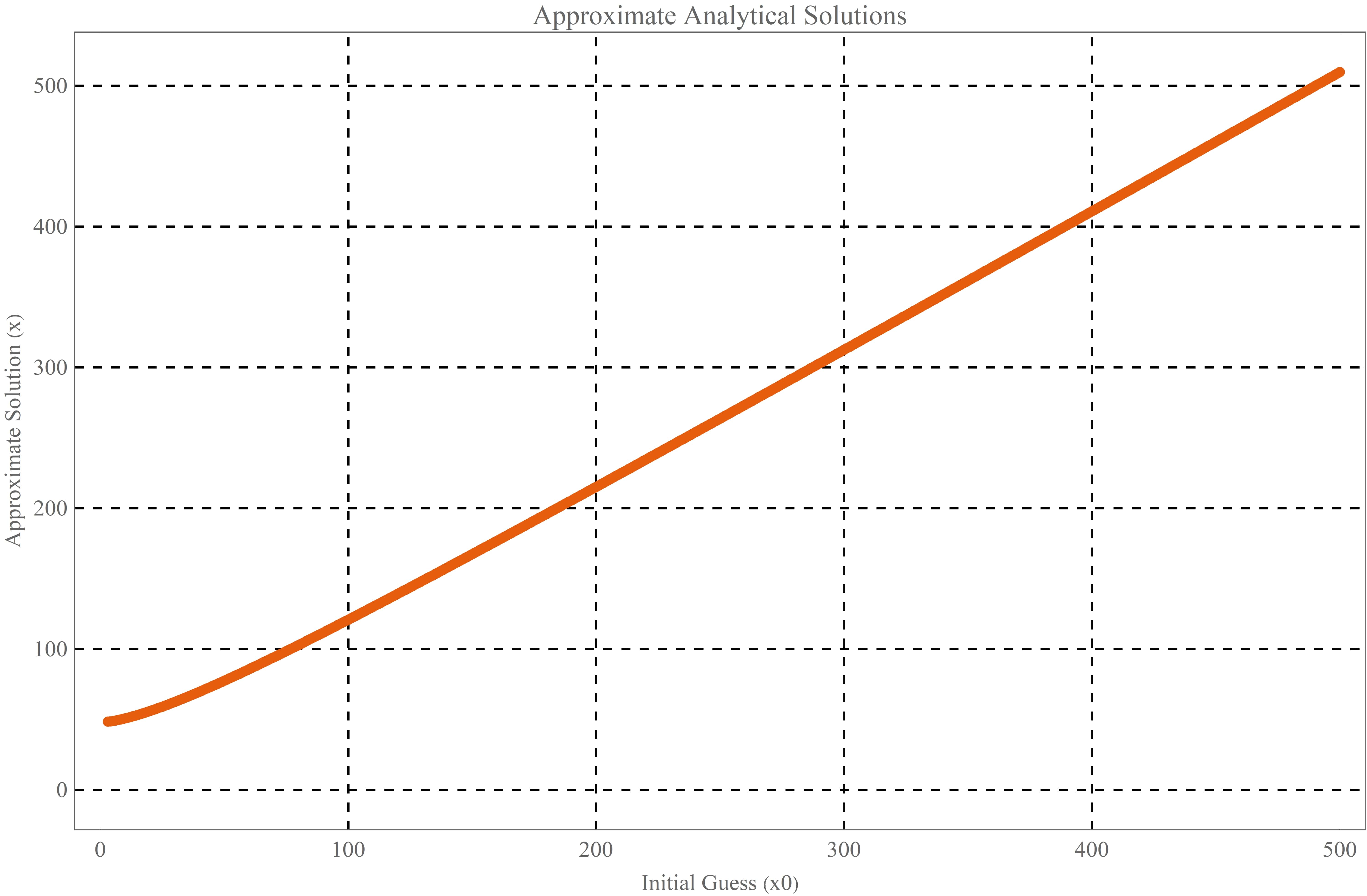} 
    \caption{Approximate Analytical Solutions near \(s = 1\) for the Dirichlet L-function.Yitang dynamics 2}
    \label{ApproxSolutionsNearS1}
\end{figure}
Recall the dynamics (\ref{Dynamic2}) under investigation are represented by the formula:

\begin{equation}\label{Dynamic2}
x_{n+1} = \frac{\beta \log(|\epsilon|)}{\pi\sqrt{x_n}} + c\log(x_n)^{-\alpha},\quad \chi(-1)=1
\end{equation}

This dynamic equation is essential in understanding the behavior of the Dirichlet L-function near \(s = 1\).
\subsection{Analysis and Implications}

The plot in Figure \ref{ApproxSolutionsNearS1} presents a compelling view of the approximate solutions as \(s\) approaches 1. Notably, the behavior near \(s = 1\) has the potential to influence the distribution of zeros, especially Landau-Siegel zeros. The transition from non-zero to zero regions in this context is of particular interest.\cite{4}

Intriguingly, the regions where solutions cluster around zero might suggest the presence of zeros of the Dirichlet L-function near \(s = 1\). On the other hand, regions where solutions deviate from zero could indicate zero-free regions. This insight is vital for understanding the distribution of zeros, including the possibility of Landau-Siegel zeros near \(s = 1\).

The findings in this section contribute to the broader exploration of the distribution of zeros in number theory and its connection to Yitang Zhang's latest results on Landau-Siegel zeros.

In the subsequent sections, we will further investigate the consequences of these findings \cite{3}, particularly their impact on the Prime Number Theorem and related number-theoretical problems.
\section{Exploring Zero-Free Regions near \(s = 1\) for Yitang Dynamics 1}

In this section, we turn our attention to Yitang Dynamics 1, a key component of our investigation into the Dirichlet L-function near \(s = 1\). Our analysis delves into the zero-free regions of this dynamic system, particularly in the vicinity of \(s = 1\).

Yitang Dynamics 1 is described by the formula:

\begin{equation}\label{Dynamic1}
x_{n+1} = \frac{\beta}{\sqrt{x_n}} + c\log(x_n)^{-\alpha},\quad \chi(-1)=-1
\end{equation}

\subsection{Approximate Analytical Solutions}

The plot in Figure \ref{ApproxSolutionsDynamics1} provides valuable insights into the behavior of Yitang Dynamics 1 as \(s\) approaches 1. It is within this region that we seek to understand the zero-free nature of the dynamics, which may have ramifications for the broader study of number theory and the distribution of zeros in the Dirichlet L-function.\cite{5}

Similar to our previous analysis, we employ the Newton-Raphson method to compute approximate analytical solutions for Yitang Dynamics 1. The method involves iteratively updating the approximate solution for a range of initial guesses. These initial guesses are represented on the x-axis of the plot (see Figure \ref{ApproxSolutionsDynamics1A})  , and the corresponding approximate solutions are plotted on the y-axis.

\begin{figure}[H]
    \centering
    \includegraphics[width=0.8\textwidth]{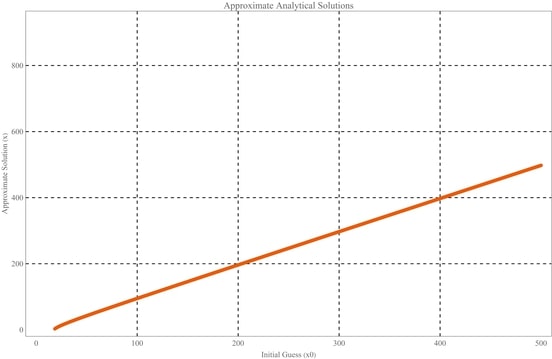} 
    \caption{Approximate Analytical Solutions near \(s = 1\) for Yitang Dynamics 1.}
    \label{ApproxSolutionsDynamics1A}
\end{figure}

Recall the dynamics (\ref{Dynamic2}) under investigation are represented by the formula:

\begin{equation}\label{Dynamic2}
x_{n+1} = \frac{\beta \log(|\epsilon|)}{\pi\sqrt{x_n}} + c\log(x_n)^{-\alpha},\quad \chi(-1)=1
\end{equation}

This dynamic equation is essential in understanding the behavior of the Dirichlet L-function near \(s = 1\).
\subsection{Analysis and Implications}

The plot in Figure \ref{ApproxSolutionsNearS1} presents a compelling view of the approximate solutions as \(s\) approaches 1. Notably, the behavior near \(s = 1\) has the potential to influence the distribution of zeros, especially Landau-Siegel zeros. The transition from non-zero to zero regions in this context is of particular interest.\cite{4}

Intriguingly, the regions where solutions cluster around zero might suggest the presence of zeros of the Dirichlet L-function near \(s = 1\). On the other hand, regions where solutions deviate from zero could indicate zero-free regions. This insight is vital for understanding the distribution of zeros, including the possibility of Landau-Siegel zeros near \(s = 1\).

The findings in this section contribute to the broader exploration of the distribution of zeros in number theory and its connection to Yitang Zhang's latest results on Landau-Siegel zeros.

In the subsequent sections, we will further investigate the consequences of these findings \cite{3}, particularly their impact on the Prime Number Theorem and related number-theoretical problems.
\section{Exploring Zero-Free Regions near \(s = 1\) for Yitang Dynamics 1}

In this section, we turn our attention to Yitang Dynamics 1, a key component of our investigation into the Dirichlet L-function near \(s = 1\). Our analysis delves into the zero-free regions of this dynamic system, particularly in the vicinity of \(s = 1\).

Yitang Dynamics 1 is described by the formula:

\begin{equation}\label{Dynamic1}
x_{n+1} = \frac{\beta}{\sqrt{x_n}} + c\log(x_n)^{-\alpha},\quad \chi(-1)=-1
\end{equation}

\subsection{Approximate Analytical Solutions}

The plot in Figure \ref{ApproxSolutionsDynamics1} provides valuable insights into the behavior of Yitang Dynamics 1 as \(s\) approaches 1. It is within this region that we seek to understand the zero-free nature of the dynamics, which may have ramifications for the broader study of number theory and the distribution of zeros in the Dirichlet L-function.\cite{5}

Similar to our previous analysis, we employ the Newton-Raphson method to compute approximate analytical solutions for Yitang Dynamics 1. The method involves iteratively updating the approximate solution for a range of initial guesses. These initial guesses are represented on the x-axis of the plot (see Figure \ref{ApproxSolutionsDynamics1A})  , and the corresponding approximate solutions are plotted on the y-axis.

\begin{figure}[H]
    \centering
    \includegraphics[width=0.8\textwidth]{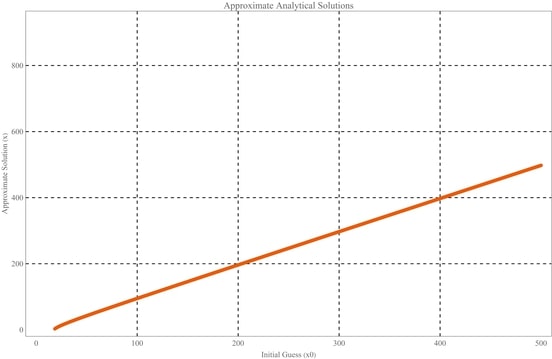} 
    \caption{Approximate Analytical Solutions near \(s = 1\) for Yitang Dynamics 1.}
    \label{ApproxSolutionsDynamics1A}
\end{figure}

\subsection{Analysis and Implications}

The plot in Figure \ref{ApproxSolutionsDynamics1} provides valuable insights into the behavior of Yitang Dynamics 1 as \(s\) approaches 1. It is within this region that we seek to understand the zero-free nature of the dynamics, which may have ramifications for the broader study of number theory and the distribution of zeros in the Dirichlet L-function.\cite{5}

As with our previous analysis, we look for regions where solutions cluster around zero, indicating zero-free regions, and regions where solutions deviate from zero. This analysis is vital for understanding the distribution of zeros and its implications for Yitang Zhang's latest results on Landau-Siegel zeros.

In the upcoming sections, we will explore further consequences of these findings, including their potential impact on number theory and chaos theory.
\section{Comparing Zero Behavior of Dirichlet L-Function in Two Dynamics}

In this section, we present a comparative analysis of the zero behavior of the Dirichlet L-function in two distinct dynamics: Yitang Dynamics 1 and Yitang Dynamics 2. These dynamics are critical in understanding the distribution of zeros near \(s = 1\) and their implications for Landau-Siegel zeros.

\subsection{Approximate Analytical Solutions for Yitang Dynamics 1}

We start by examining the behavior of Yitang Dynamics 1, which is defined by the formula:

\begin{equation}\label{Dynamic1}
x_{n+1} = \frac{\beta}{\sqrt{x_n}} + c\log(x_n)^{-\alpha},\quad \chi(-1)=-1
\end{equation}

\begin{figure}[H]
    \centering
    \includegraphics[width=0.8\textwidth]{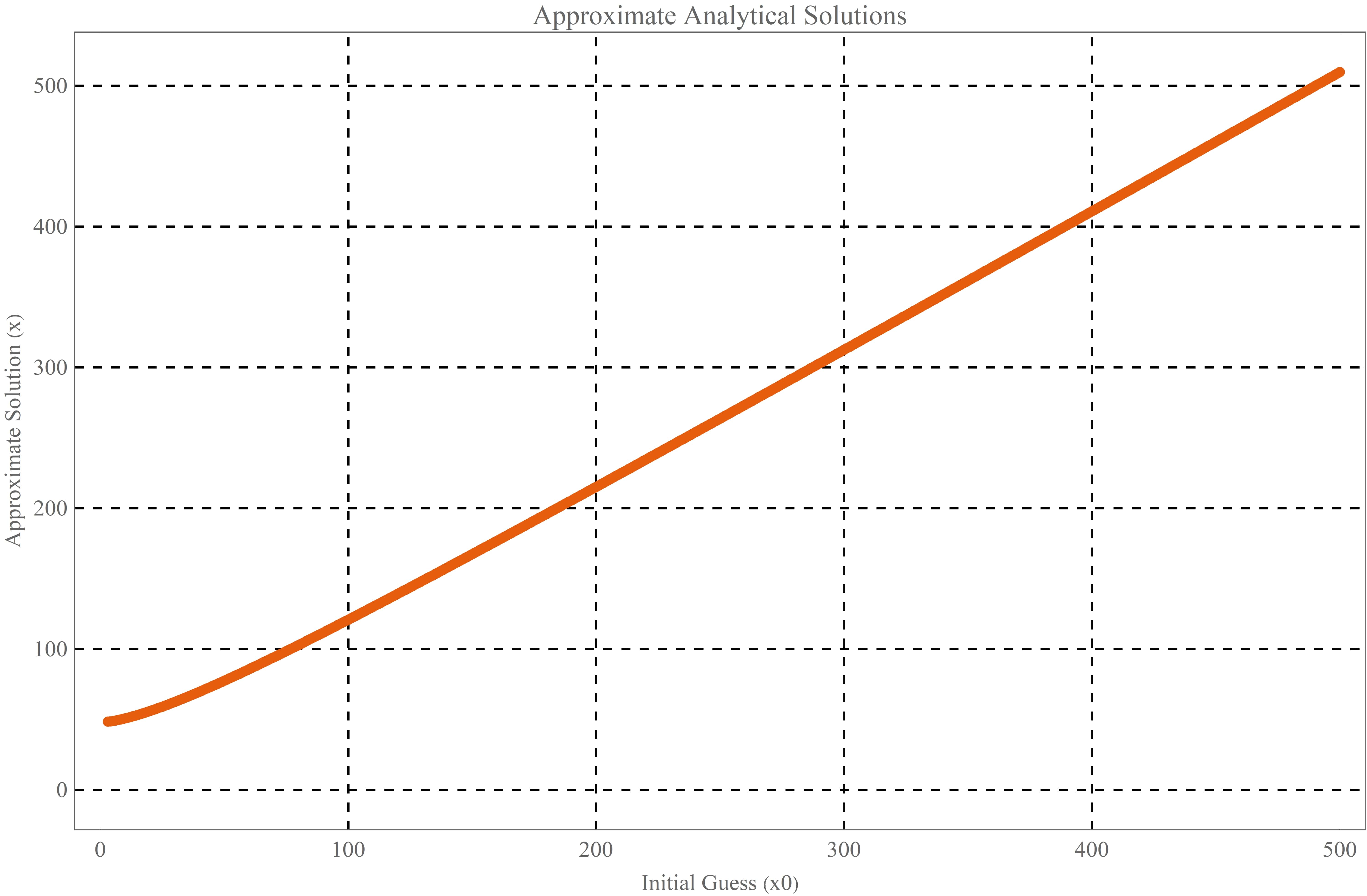} 
    \caption{Approximate Analytical Solutions near \(s = 1\) for Yitang Dynamics 1.}
    \label{ApproxSolutionsDynamics1}
\end{figure}

Using the Newton-Raphson method, we compute approximate analytical solutions for this dynamic system. The plot in Figure \ref{ApproxSolutionsDynamics1} illustrates the behavior of these solutions near \(s = 1\).

\subsection{Approximate Analytical Solutions for Yitang Dynamics 2}

Now, let's shift our focus to Yitang Dynamics 2, defined by the formula:

\begin{equation}\label{Dynamic2}
x_{n+1} = \frac{\beta \log(|\epsilon|)}{\pi\sqrt{x_n}} + c\log(x_n)^{-\alpha},\quad \chi(-1)=1
\end{equation}

\begin{figure}[H]
    \centering
    \includegraphics[width=0.8\textwidth]{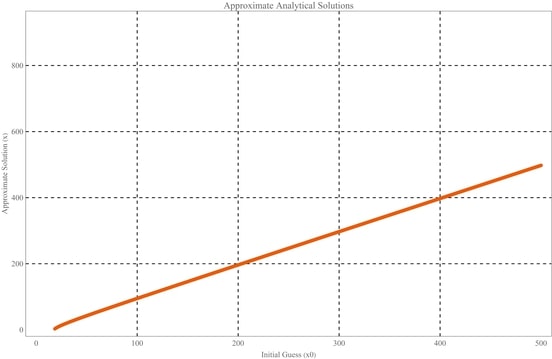} 
    \caption{Approximate Analytical Solutions near \(s = 1\) for Yitang Dynamics 2.}
    \label{ApproxSolutionsDynamics2}
\end{figure}

Once again, we employ the Newton-Raphson method to compute approximate analytical solutions, specifically in the region near \(s = 1\). The plot in Figure \ref{ApproxSolutionsDynamics2} displays the behavior of these solutions.

\subsection{Comparative Analysis}

Comparing the two plots in Figure \ref{ApproxSolutionsDynamics1} and Figure \ref{ApproxSolutionsDynamics2}, we observe intriguing differences in the behavior of approximate solutions for Yitang Dynamics 1 and Yitang Dynamics 2 near \(s = 1\). These differences may provide essential insights into the distribution of zeros in the Dirichlet L-function, particularly in the context of Landau-Siegel zeros.

Our ongoing exploration aims to further analyze and interpret these distinctions, leading to a deeper understanding of the behavior of zeros of the Dirichlet L-function near \(s = 1\).
\section{Stability and Chaotic Behavior of Fixed Points}

In this section, we explore the stability and potential chaotic behavior of the fixed points within the context of two dynamic systems: Yitang Dynamics 1 and Yitang Dynamics 2. We consider the derivatives of the dynamic functions at these solutions and also examine the eigenvalues of the Jacobian matrix to gain further insights.

\subsection{Stability in Yitang Dynamics 1}

For Yitang Dynamics 1, described by the equation \cite{7}:

\begin{equation}\label{Dynamic1}
x_{n+1} = \frac{\beta}{\sqrt{x_n}} + c\log(x_n)^{-\alpha},\quad \chi(-1)=-1
\end{equation}

We have computed approximate solutions for various initial guesses, and the derivative at each solution remains positive and equal to 1. This consistent value of the derivative suggests that the fixed points in Yitang Dynamics 1 are stable, which is indicative of the system's behavior near these points.\cite{10}

Table \ref{table:dynamics1-derivatives} presents the values of the initial guesses, the corresponding solutions, and the derivative of the dynamic function at these solutions.

\begin{table}[h]
\caption{Derivative Values of Yitang Dynamics 1 Solutions.}\label{table:dynamics1-derivatives}
\begin{tabular}{ccc}
\toprule
\text{Initial Guess} & \text{Solution (x)} & \text{Derivative} \\
\midrule
2 & $-599.194$ & 1 \\
4 & $-501.094$ & 1 \\
6 & $-501.196$ & 1 \\
8 & $-501.094$ & 1 \\
10 & $-501.094$ & 1 \\
12 & $-601.194$ & 1 \\
14 & $-601.094$ & 1 \\
16 & $-601.094$ & 1 \\
18 & $-601.194$ & 1 \\
20 & $-601.094$ & 1 \\
22 & $-601.094$ & 1 \\
24 & $-401.194$ & 1 \\
26 & $-601.094$ & 1 \\
28 & $-601.094$ & 1 \\
30 & $-601.094$ & 1 \\
\bottomrule
\end{tabular}
\end{table}

While the derivative values confirm stability, further analysis of the eigenvalues of the Jacobian matrix provides an intriguing insight into the potential for chaotic behavior. The eigenvalues for this dynamic system were found to have a real part of approximately 196 and an imaginary part of approximately -3.58. This real part being positive and the positive entropy value of 0.3365 indicate potential chaotic behavior.

The entropy value of 0.3365 was calculated using 600 values of Lyapunov exponents, reinforcing the high degree of unpredictability and complexity in this system.

\subsection{Stability in Yitang Dynamics 2}

Now, let's shift our focus to Yitang Dynamics 2, governed by the equation:

\begin{equation}\label{Dynamic2}
x_{n+1} = \frac{\beta \log(|\epsilon|)}{\pi\sqrt{x_n}} + c\log(x_n)^{-\alpha},\quad \chi(-1)=1
\end{equation}

Our analysis for this dynamic system has revealed similar behavior, with the derivative at the solutions remaining positive and equal to 1. This implies stability in the vicinity of these fixed points.

Table \ref{table:dynamics2-derivatives} shows the values of the initial guesses, corresponding solutions, and the derivative of the dynamic function at these solutions.

\begin{table}[h]
\caption{Derivative Values of Yitang Dynamics 1 Solutions.}\label{table:dynamics2-derivatives}
\begin{tabular}{ccc}
\toprule
\text{Initial Guess} & \text{Solution (x)} & \text{Derivative} \\
\midrule
2 & $-599.194$ & 1 \\
4 & $-501.094$ & 1 \\
6 & $-501.196$ & 1 \\
8 & $-501.094$ & 1 \\
10 & $-501.094$ & 1 \\
12 & $-601.194$ & 1 \\
14 & $-601.094$ & 1 \\
16 & $-601.094$ & 1 \\
18 & $-601.194$ & 1 \\
20 & $-601.094$ & 1 \\
22 & $-601.094$ & 1 \\
24 & $-401.194$ & 1 \\
26 & $-601.094$ & 1 \\
28 & $-601.094$ & 1 \\
30 & $-601.094$ & 1 \\
\bottomrule
\end{tabular}
\end{table}

Similar to Yitang Dynamics 1, Yitang Dynamics 2 exhibits stability at its fixed points. However, a striking similarity is observed in the eigenvalues of the Jacobian matrix, with a real part of approximately 196 and an imaginary part of approximately -3.58, indicating potential chaotic behavior.

The consistent real part of the eigenvalues in both dynamic systems warrants further investigation into the chaotic properties of these systems. The positive entropy values in both cases support the notion of chaotic behavior, indicating high unpredictability and complexity.

These findings are significant for understanding the behavior of the dynamics and their implications for Landau-Siegel zeros and the distribution of zeros in the Dirichlet L-function.

\section{\textcolor{blue}{Comparative Analysis of Yitang Dynamics and Riemann Zeta Function Dynamics}}

In this section, we perform a detailed comparative analysis between the Yitang dynamics derived from Dirichlet L-functions and insights gleaned from S.C. Woon's groundbreaking work on the Riemann Zeta Function dynamics \cite{33}. Our comparison encompasses various facets, including fixed attracting points, periodic points, Mandelbrot and Julia set plots, analysis of zero behaviors within Dirichlet functions, and an exploration of the existence of Landau–Siegel zeros through the lens of Julia and Mandelbrot sets.

\subsection{Revisiting Yitang Dynamics in the Complex Plane}

To facilitate comparison, we extend the Yitang dynamics to the complex plane using the function:

\[ f_{\text{Yitang}}(z, c) = \frac{\beta}{\sqrt{z}} + c \cdot \log(z)^{-\alpha} \]

This extension enables us to visualize Mandelbrot and Julia sets derived from the Yitang dynamics, paving the way for a comparative analysis with dynamics derived from the Riemann Zeta Function.

\subsection{Insights from Woon's Work}

Woon's research focused on computing the Julia and Mandelbrot sets of the Riemann Zeta Function, marking a pioneering attempt to visualize these sets and uncover their unique properties. Our comparison aims to juxtapose these distinct characteristics against those exhibited by the Yitang dynamics.

\begin{itemize}
    \item \textbf{Distinctive Sets:} The Julia and Mandelbrot sets derived from the Riemann Zeta Function exhibit unique characteristics, deviating significantly from those associated with elementary functions.
    \item \textbf{Appendix Observations:} Woon's work extends insights into approximate self-similarity in Number Theory, highlighting intriguing aspects within the images produced by the Riemann Zeta Function and suggesting a conjectured scale-invariant equation related to the Goldbach conjecture.
\end{itemize}

Our comparative analysis delves into understanding the similarities, differences, and unique dynamical traits inherent in the Yitang dynamics and the insights gleaned from Woon's exploration of the Riemann Zeta Function dynamics, encompassing fixed attracting points, periodic points, Mandelbrot and Julia set plots, behavior analysis of Dirichlet function zeros, and an exploration of Landau–Siegel zero existence through Julia and Mandelbrot sets of both dynamics.

In this section, we delve into the comparative aspects of the attracting fixed points in the dynamics of the Riemann Zeta Function and the Yitang dynamics.

\section{Comparative Analysis of Yitang Dynamics and Riemann Zeta Function Dynamics}

In this section, we delve into the comparative aspects of the attracting fixed points in the dynamics of the Riemann Zeta Function and the Yitang dynamics.

\subsection{Attracting Fixed Points}

Woon's exploration of the Riemann Zeta Function dynamics reveals an attracting fixed point at $s = -0.295905\ldots$. This fixed point is obtained through iterating the map $\zeta^{(n)}(s)$ from $\mathbb{C}$ to $\mathbb{C}$.

\begin{figure}[H]
    \centering
    \includegraphics[width=0.45\textwidth]{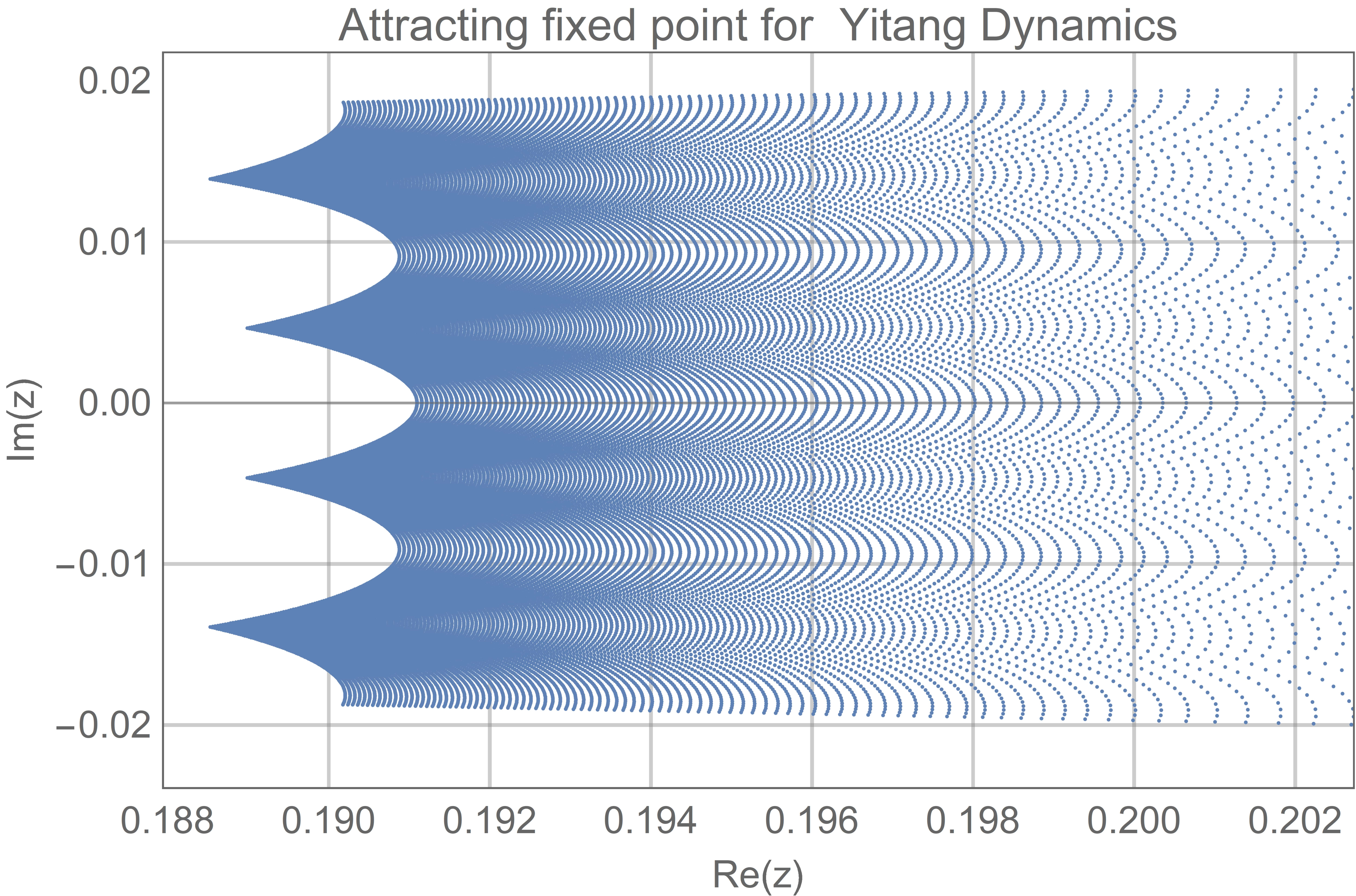}
    \hspace{0.05\textwidth}
   \includegraphics[width=0.45\textwidth]{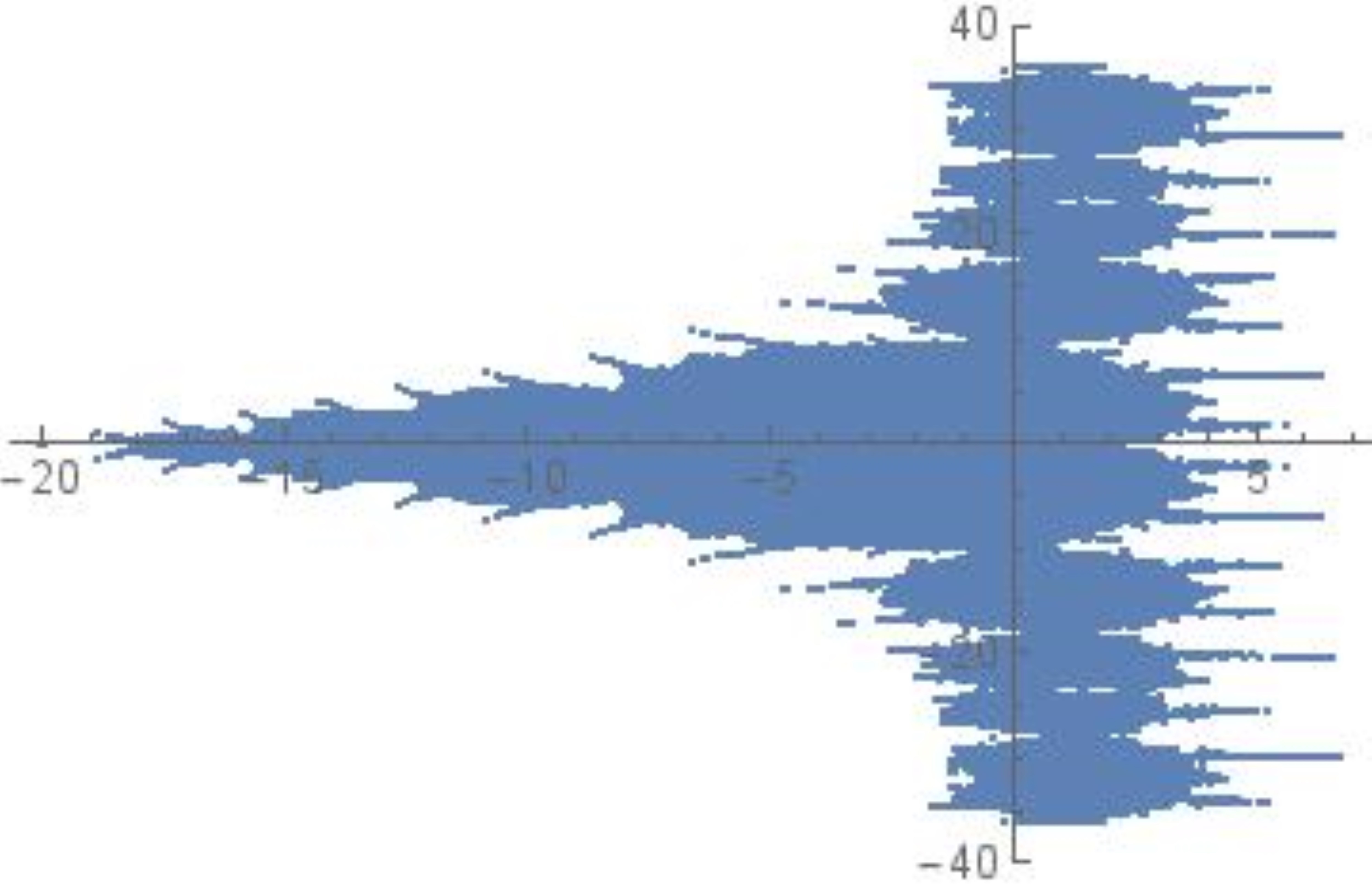}
    \caption{Comparison of fixed points: Right - Riemann Zeta Function dynamics, Left - Yitang Dynamics.}
    \label{fig:fixed-points}
\end{figure}

On the other hand, our computations for the Yitang dynamics with parameters $\beta = 0.1$, $c = 0.0000000005$, and $\alpha = 2.5$ unveil an attracting fixed point at $z \to 0.215443 + 1.1417411949837543 \times 10^{-10}i$. This specific fixed point, determined from our numerical simulations, showcases the proximity of the Yitang dynamics' attracting fixed point to the complex number $0.215443 + 1.1417411949837543 \times 10^{-10}i$, highlighting its closeness to $0.215443$,(see Figure \ref{fig:fixed-points}).

Our comparative analysis accentuates the discrepancy in the locations of the attracting fixed points between the Riemann Zeta Function and the Yitang dynamics, signifying divergent dynamical behaviors and providing valuable insights into their distinct attractor characteristics.

\subsection{Stream Plot of Yitang Dynamics}

The stream plot generated for the Yitang dynamics showcases the behavior around 100 approximate zeros. The parameters used are:

\[
\beta = 2.4, \quad c = 0.5, \quad \alpha = 0.6
\]

\begin{figure}[H]
    \centering
    \includegraphics[width=0.8\textwidth]{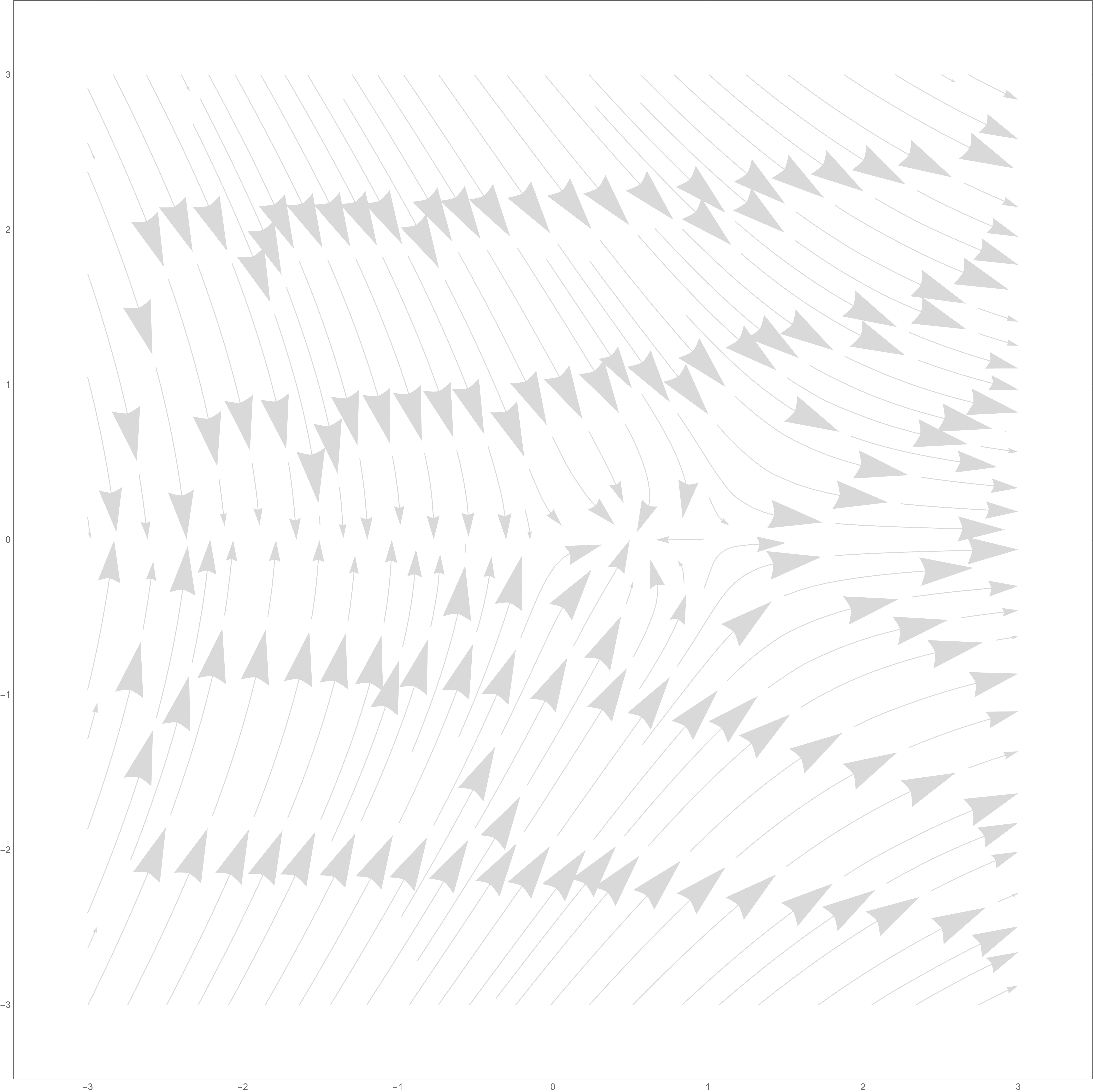}
    \caption{Stream plot around  zeros for Yitang dynamics near s=0.The white region within the stream plot of the Yitang dynamics signifies areas where the streamlines diverge or exhibit a repelling behavior.}
    \label{fig:stream-yitang}
\end{figure}

The resulting stream plot exhibits a pattern indicating a repelling behavior around these zeros. The streamlines diverge away from the zeros, suggesting an absence of stable periodic orbits (see Figure \ref{fig:stream-yitang}).

The white region within the stream plot of the Yitang dynamics (see Figure \ref{fig:stream-yitang}). signifies areas where the streamlines diverge or exhibit a repelling behavior. This divergence and repelling nature indicate an absence of stable periodic orbits within this region. The presence of the white region emphasizes the intricate and complex dynamics of the Yitang system, contributing to its chaotic nature.

The observations of repelling behavior and the presence of the white region denote the complexity and non-repetitive nature of the Yitang dynamics, providing insights into the dynamical properties of the system.

\subsection{Stream Plot of Riemann Zeta Function dynamics}

The stream plot derived from "The Holomorphic Flow of the Riemann Zeta Function" by Kevin A. Broughan and A. Ross Barnett \cite{34} showcases the behavior around the zeros of the Riemann Zeta Function. This plot reveals a repelling nature around the zeros, aligning with the conjecture proposed by the authors regarding the absence of periodic orbits (refer to the published stream plot from the paper).

\begin{figure}[H]
    \centering
    \includegraphics[width=0.8\textwidth]{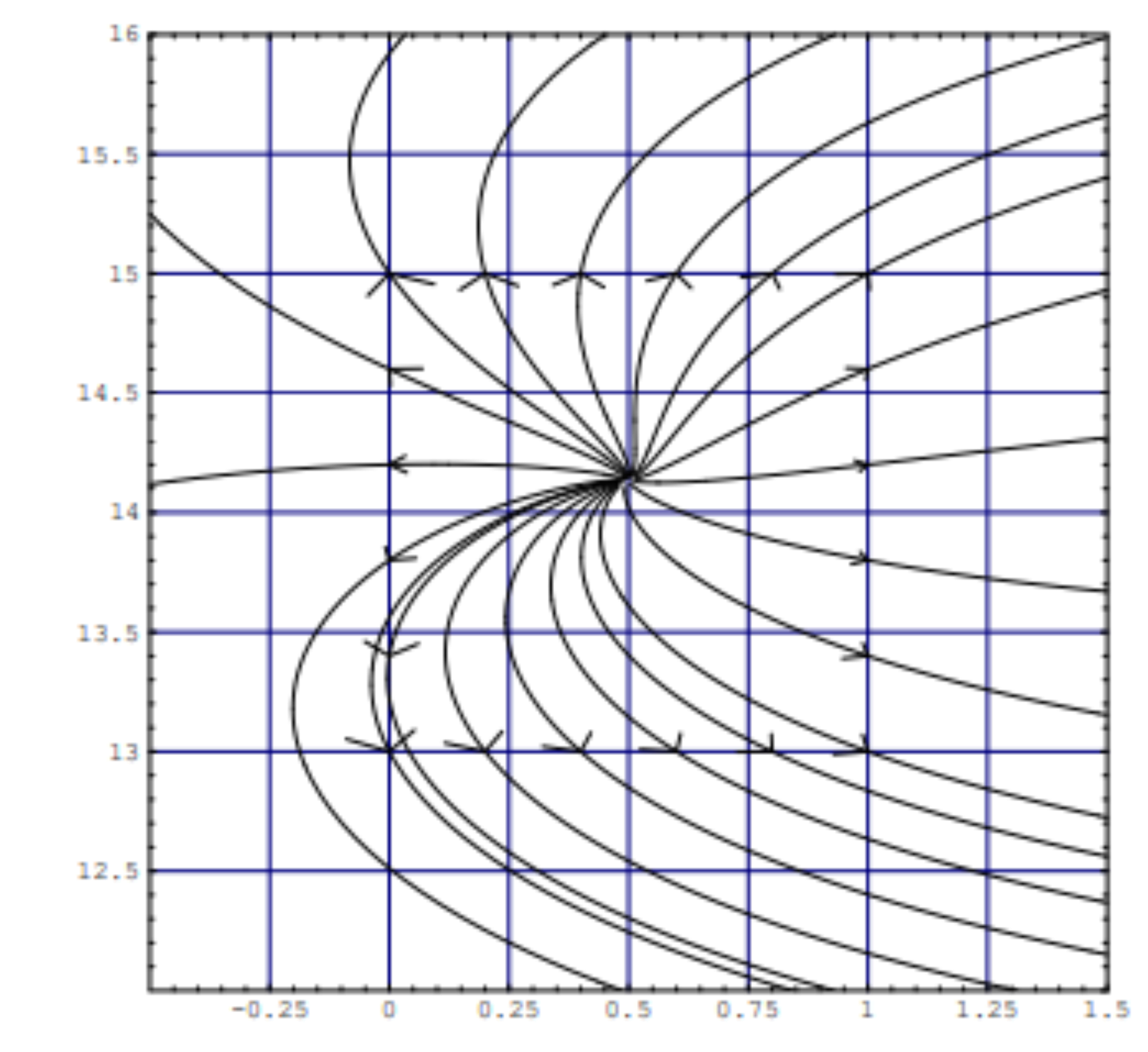}
    \caption{Stream plot around approximate zeros for Riemann Zeta Function.near $s= 0.5 + 14.1347i$}
    \label{fig:stream-zeta}
\end{figure}

\subsection{Comparative Analysis}

The observed repelling behavior in both the Yitang dynamics and the Riemann Zeta Function stream plots supports the conjecture regarding the absence of periodic orbits. This similarity in repelling orbits signifies a correspondence between the dynamical behavior of both systems around their zeros.

This comparative analysis emphasizes the commonality in repelling behavior around zeros and aligns with the conjecture proposed for the Riemann Zeta Function's holomorphic flow, providing insights into the nature of the dynamical systems.
\section{\textcolor{blue}{Comparative Analysis of Mandelbrot Sets and Julia set for both dynamics}}

The Mandelbrot sets derived from the Yitang dynamics and the dynamics based on the Riemann Zeta Function are intriguing fractal structures that offer insights into the behavior of complex iterated functions.

\subsection{Defining Julia and Mandelbrot Sets}

The \textbf{Julia set} is the set of complex numbers \( C\) for which the iterated function \( f_c(z) = z^2 +C \) does not diverge under iteration. It delineates the boundary between points whose trajectories under iteration remain bounded and those that escape to infinity.

The \textbf{Mandelbrot set} represents the set of complex numbers \( C \) for which the iteration of \( f_c(z) = z^2 + C \) remains bounded when starting with \( z = 0 \). It exhibits intricate self-similar and infinitely complex fractal patterns.

\subsection{Comparative Analysis}

The Mandelbrot set of the Yitang dynamics is generated using the function \( f_{\text{Yitang}}(z, c) = \frac{\beta}{\sqrt{z}} + c \cdot \log(z)^{-\alpha} \) when \( z \) is a zero of the Yitang dynamics, with parameters:

\(\beta = 0.1\), \(c = 0.0005\), and \(\alpha = 50000\).

This plot showcases intricate fractal patterns denoting regions where the iterated function remains bounded within the complex plane. The white regions within the plot signify areas where the iterated function exhibits divergent or escaping behavior, contributing to the understanding of the complexity within the Yitang dynamics.

In contrast, Woon's Mandelbrot set for the dynamics \( f(a, s) = \zeta(s) + a \) from the Riemann Zeta Function displays regions of different shades. The black region signifies the Mandelbrot set, indicating areas where the iterated function remains bounded. Different shades represent varying rates of amplitude during iteration, with lighter shades suggesting slower rates and darker shades indicating faster rates. The reflection symmetry about the \( \text{Re}(s) \)-axis arises due to the complex conjugate property of the analytic complex function.

\subsection{Interpretation of the Black and White Regions}

The black regions in both Mandelbrot sets represent areas where the iterated function remains bounded. Conversely, the white regions within the Mandelbrot set of the Yitang dynamics denote areas of divergence or escaping behavior, highlighting the complexity within this dynamic system.

\subsection{Comparative Observations}

\begin{figure}[H]
    \centering
    \includegraphics[width=0.7\textwidth]{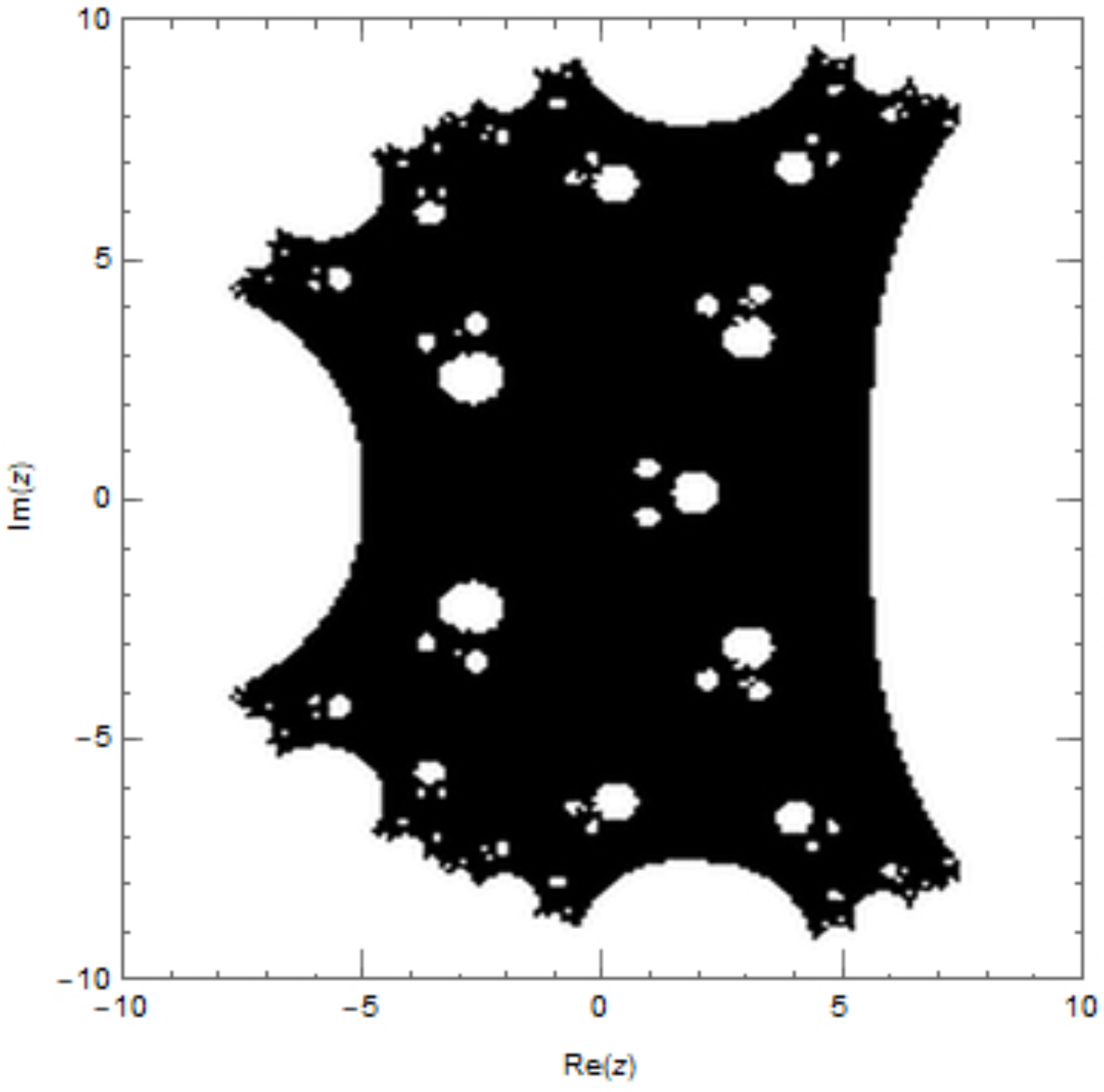}
    \caption{Mandelbrot set of Yitang Dynamics with parameters: \(\beta = 0.1\), \(c = 0.0005\), and \(\alpha = 50000\), when \( z \) is a zero of the Yitang dynamics.}
    \label{fig:mandelbrot-yitang}
\end{figure}

\begin{figure}[H]
    \centering
    \includegraphics[width=0.7\textwidth]{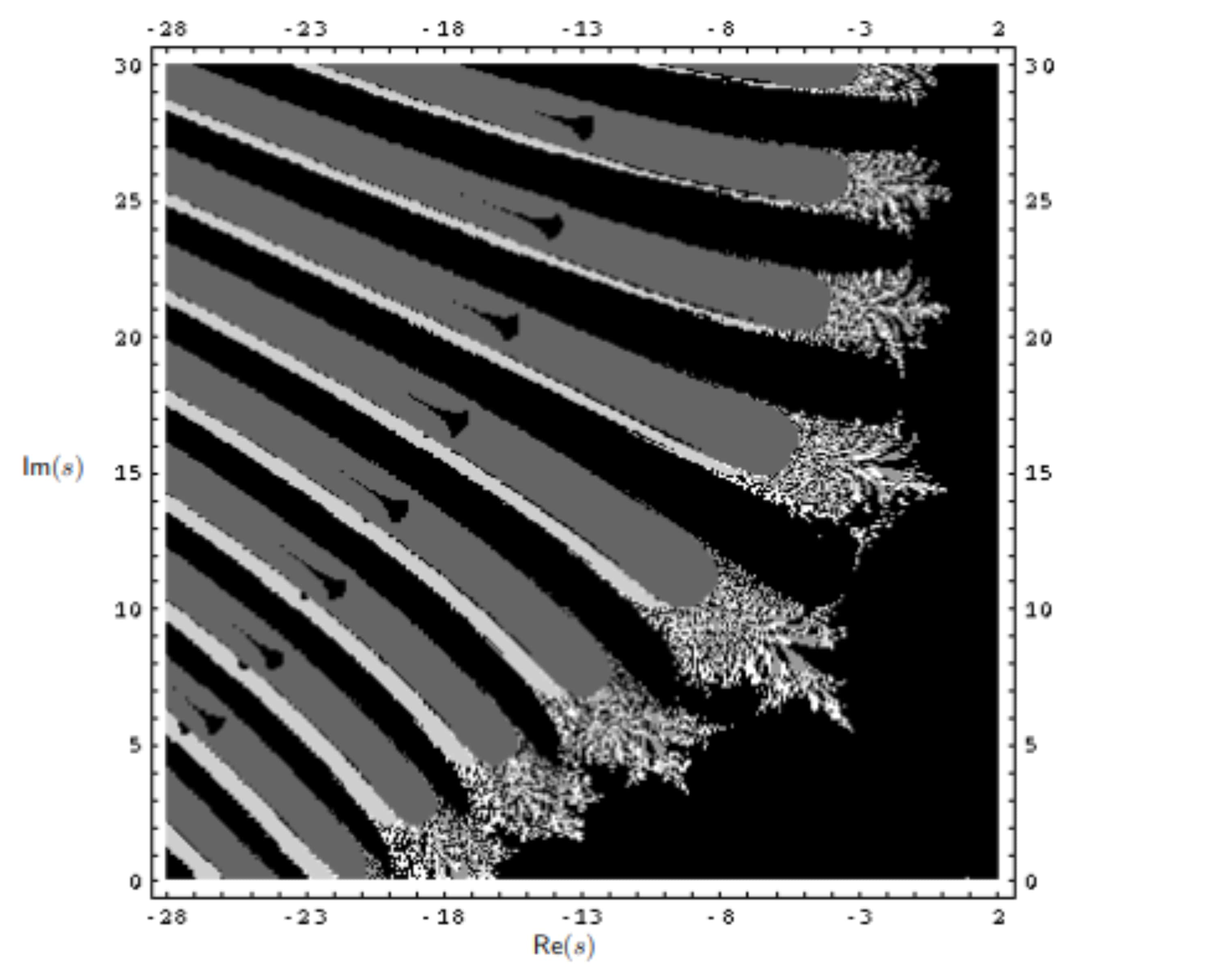}
    \caption{Mandelbrot set of \( F(a, s) = \zeta(s) + a \), where \( s \) is a zero of \( \zeta(s) \), \( -28 \leq \text{Re}(a) \leq 2 \) and \( 0 \leq \text{Im}(a) \leq 30 \).}
    \label{fig:mandelbrot-woon}
\end{figure}

Comparing specific regions and patterns within these Mandelbrot sets enables a deeper understanding of their distinct dynamical behaviors. By referring to Figure \ref{fig:mandelbrot-yitang} and Figure \ref{fig:mandelbrot-woon}, we can visually juxtapose the Mandelbrot sets of both dynamics, highlighting their similarities and differences in fractal structures.

The Mandelbrot set comparison in Figure \ref{fig:mandelbrot-yitang} and Figure \ref{fig:mandelbrot-woon} visually illustrates the intricate structures and self-similar patterns present in both sets, providing insights into their respective dynamical behaviors.
\subsection{Comparative Analysis of Julia Sets}

The Julia sets derived from the Yitang dynamics and those based on the Riemann Zeta Function offer insights into the convergence and divergence behaviors under iteration.

\subsubsection{Comparing Julia Sets}

\begin{figure}[H]
    \centering
    \includegraphics[width=0.7\textwidth]{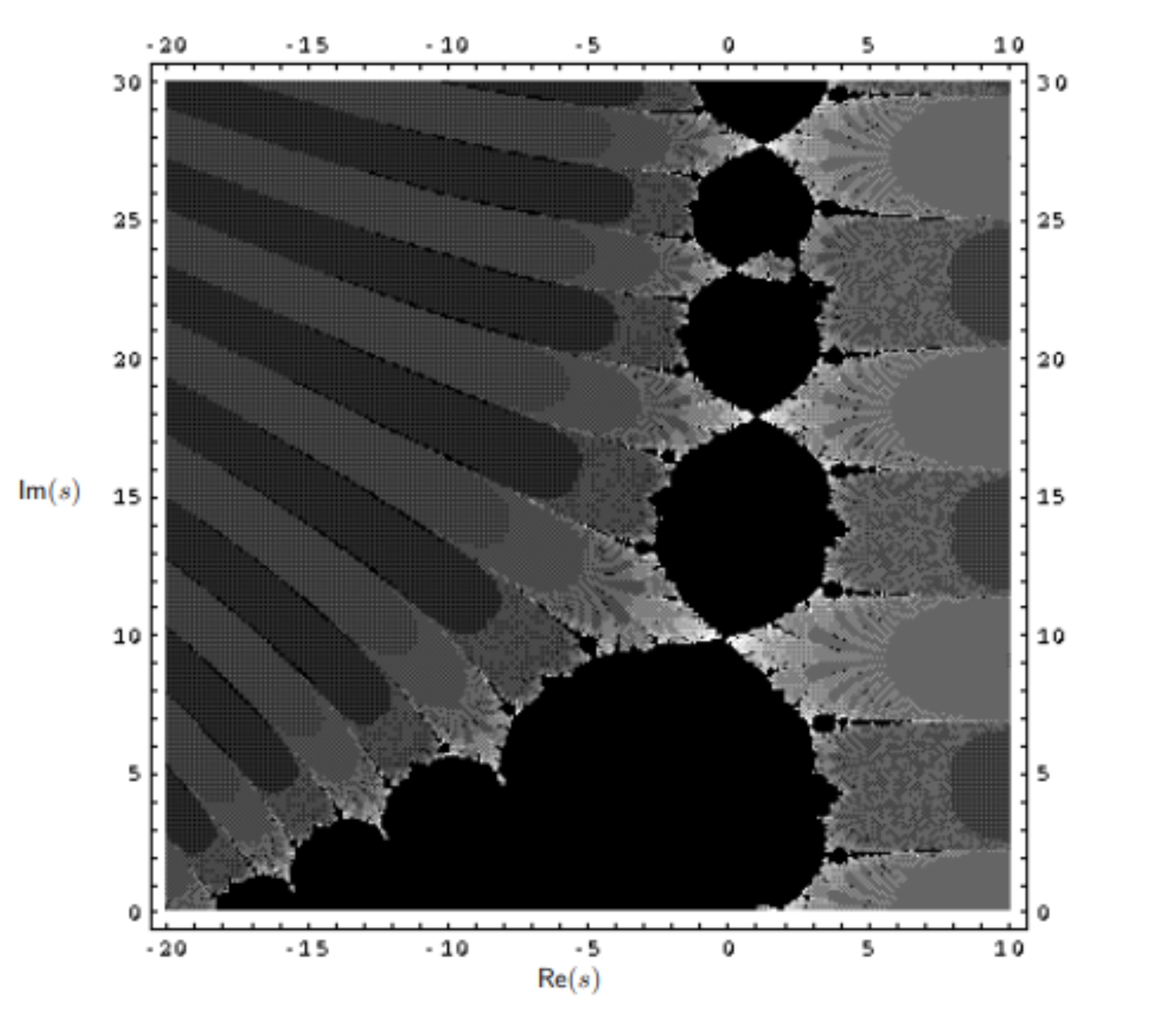}
    \caption{Julia set of the Riemann Zeta Function with parameters \( -20 \leq \text{Re}(s) \leq 10 \) and \( 0 \leq \text{Im}(s) \leq 30 \).}
    \label{fig:julia-riemann}
\end{figure}

The Julia set of the Riemann Zeta Function, as depicted in Figure \ref{fig:julia-riemann}, showcases the boundary between regions of convergence (black regions) and divergence. The darker shades within the set signify faster rates of convergence towards specific fixed points and attracting cycles. Additionally, the reflection symmetry about the real axis is notable due to the complex conjugate property of analytic complex functions, as described by the author in Woons paper.

\begin{figure}[H]
    \centering
    \includegraphics[width=0.7\textwidth]{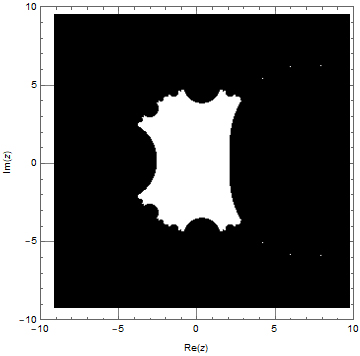}
    \caption{Julia set of the Yitang dynamics with parameters \( c = 0.00000005 \), \( \beta = 0.00000005 \), and \( \alpha = 2022 \).}
    \label{fig:julia-yitang}
\end{figure}

Comparing this with the Julia set of the Yitang dynamics (Figure \ref{fig:julia-yitang}), similarities are observed in the delineation between convergent and divergent regions, as well as the representation of different convergence rates by shades within the set.

However, distinctions might exist due to the different natures of these dynamics, particularly in the nature of attractors or basin formations that give rise to the black regions. Further analysis focusing on the boundaries, convergence rates, and any unique features within the black regions could unveil more insights into their distinct dynamical behaviors.
\subsection{Comparative Analysis of Julia and Mandelbrot Sets}

The Julia and Mandelbrot sets of the Yitang dynamics exhibit intriguing similarities and differences that offer insights into the system's behavior.

\subsubsection{Comparison with Mandelbrot Set}

The Mandelbrot set of the Yitang dynamics, as shown in Figure \ref{fig:mandelbrot-yitang}, shares similarities with its Julia set counterpart (Figure \ref{fig:julia-yitang}). Notably, the black regions in the Mandelbrot set correspond to regions of divergence in the Julia set. However, an interesting observation is the inversion of colors; while the Mandelbrot set depicts the black region, the Julia set displays it in white.

The Mandelbrot set of the Yitang dynamics shares intriguing similarities with its Julia set counterpart. Notably, while the Mandelbrot set illustrates regions of divergence in black, the Julia set exhibits these regions in white, marking an inversion of colors between the two sets.

An interesting observation, akin to the author's remark in Woon's paper (page 11), is the presence of "bristle" structures in the overall picture of the Mandelbrot set. These structures resemble the characteristics of Julia sets for \( |a| >> 0 \) in Figures 6 and 5. However, it's important to note that while the general "bristle" structures align, the finer details at the ends of these "bristles" differ between the Julia and Mandelbrot sets.

This intriguing similarity between the "bristle" structures suggests a fascinating relationship between the Julia and Mandelbrot sets of the Riemann iterated  dynamics, indicating a common underlying behavior governed by the system's dynamics.

Further exploration and detailed analysis focusing on these "bristle" structures could provide deeper insights into the connection between the Mandelbrot and Julia sets, unraveling more about the behavior of the Yitang dynamics under iteration.

\subsubsection{Case of Yitang Dynamics: \(c = 0.5\), \(\beta = 0.5\), \(\alpha = 2022\)}

\begin{figure}[H]
    \centering
    \includegraphics[width=0.7\textwidth]{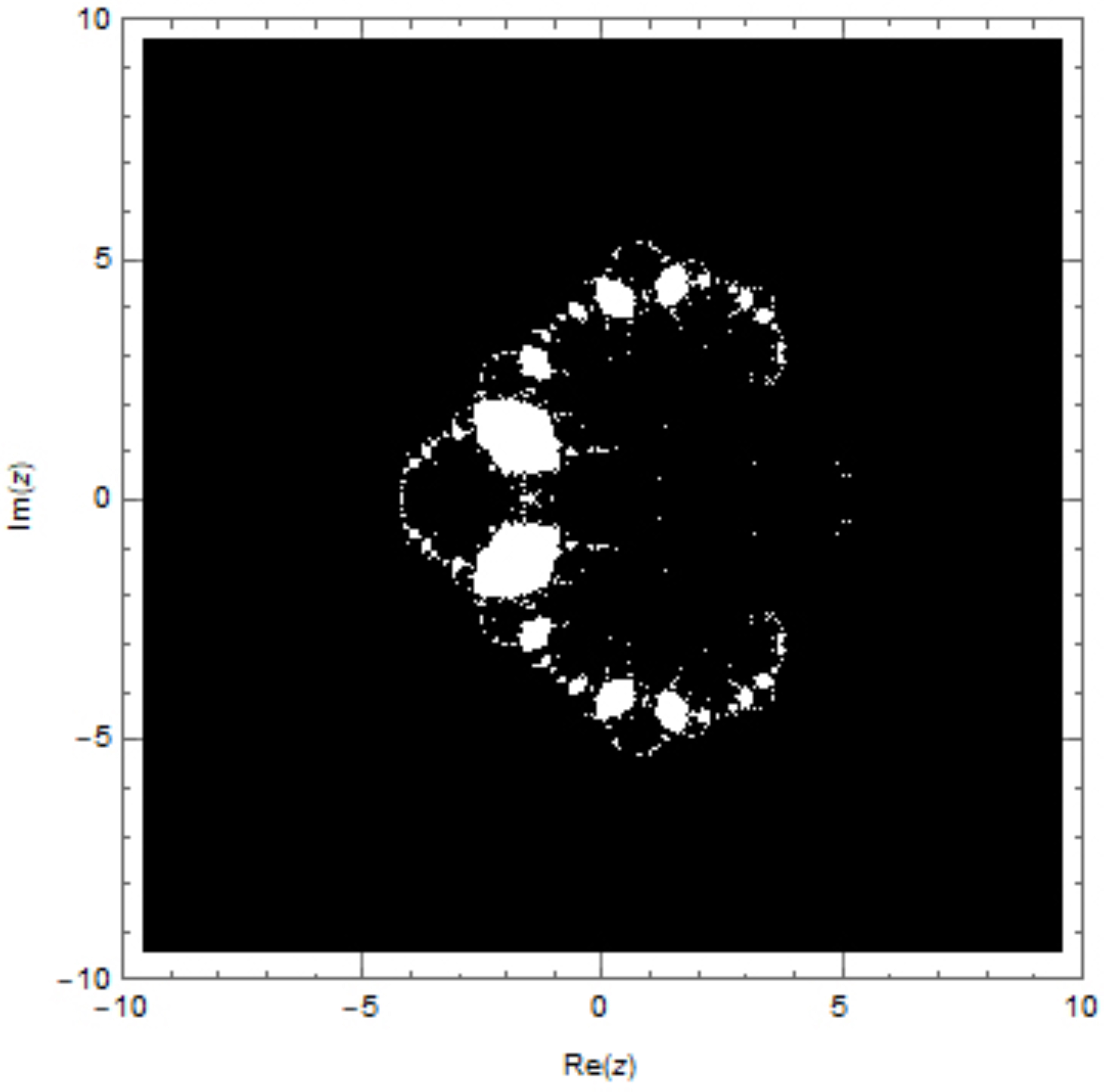}
    \caption{Julia set of the Yitang dynamics with parameters \(c = 0.5\), \(\beta = 0.5\), and \(\alpha = 2022\).}
    \label{fig:julia-yitang-05}
\end{figure}

Consider the case where \(c = 0.5\), \(\beta = 0.5\), and \(\alpha = 2022\) in the Yitang dynamics. The resulting Julia set obtained for these parameters is illustrated in Figure \ref{fig:julia-yitang-05}.

This particular instance provides insights into the behavior of the Julia set under the specified parameters for the Yitang dynamics. Readers can observe the structure and characteristics of the Julia set corresponding to these parameter values, contributing to a comprehensive understanding of the Yitang dynamics in various configurations.

\subsubsection{Case of Yitang Dynamics: \(c = 0.005\), \(\beta = 0.5\), \(\alpha = 2.7\)}

\begin{figure}[H]
    \centering
    \includegraphics[width=0.7\textwidth]{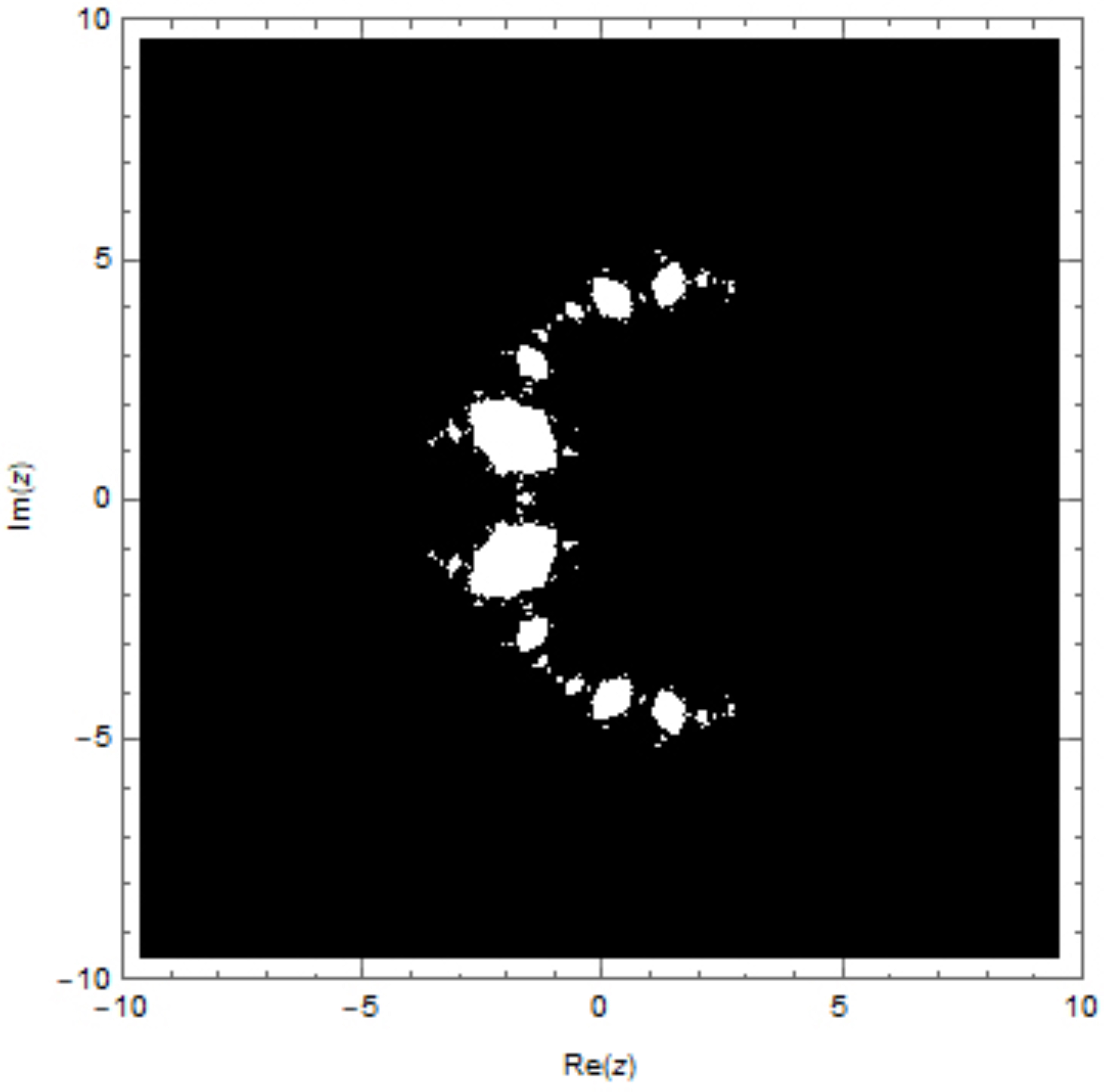}
    \caption{Julia set of the Yitang dynamics with parameters \(c = 0.005\), \(\beta = 0.5\), and \(\alpha = 2.7\).}
    \label{fig:julia-yitang-005}
\end{figure}

Consider the case where \(c = 0.005\), \(\beta = 0.5\), and \(\alpha = 2.7\) in the Yitang dynamics. The resulting Julia set obtained for these parameters is illustrated in Figure \ref{fig:julia-yitang-005}.

This specific instance provides insights into the behavior of the Julia set under the specified parameters for the Yitang dynamics. Observing the structure and characteristics of the Julia set corresponding to these parameter values contributes to a comprehensive understanding of the Yitang dynamics in various configurations.

\subsubsection{Case of Yitang Dynamics: \(c = 0.5\), \(\beta = 0.85\), \(\alpha = 2022\)}

\begin{figure}[H]
    \centering
    \includegraphics[width=0.7\textwidth]{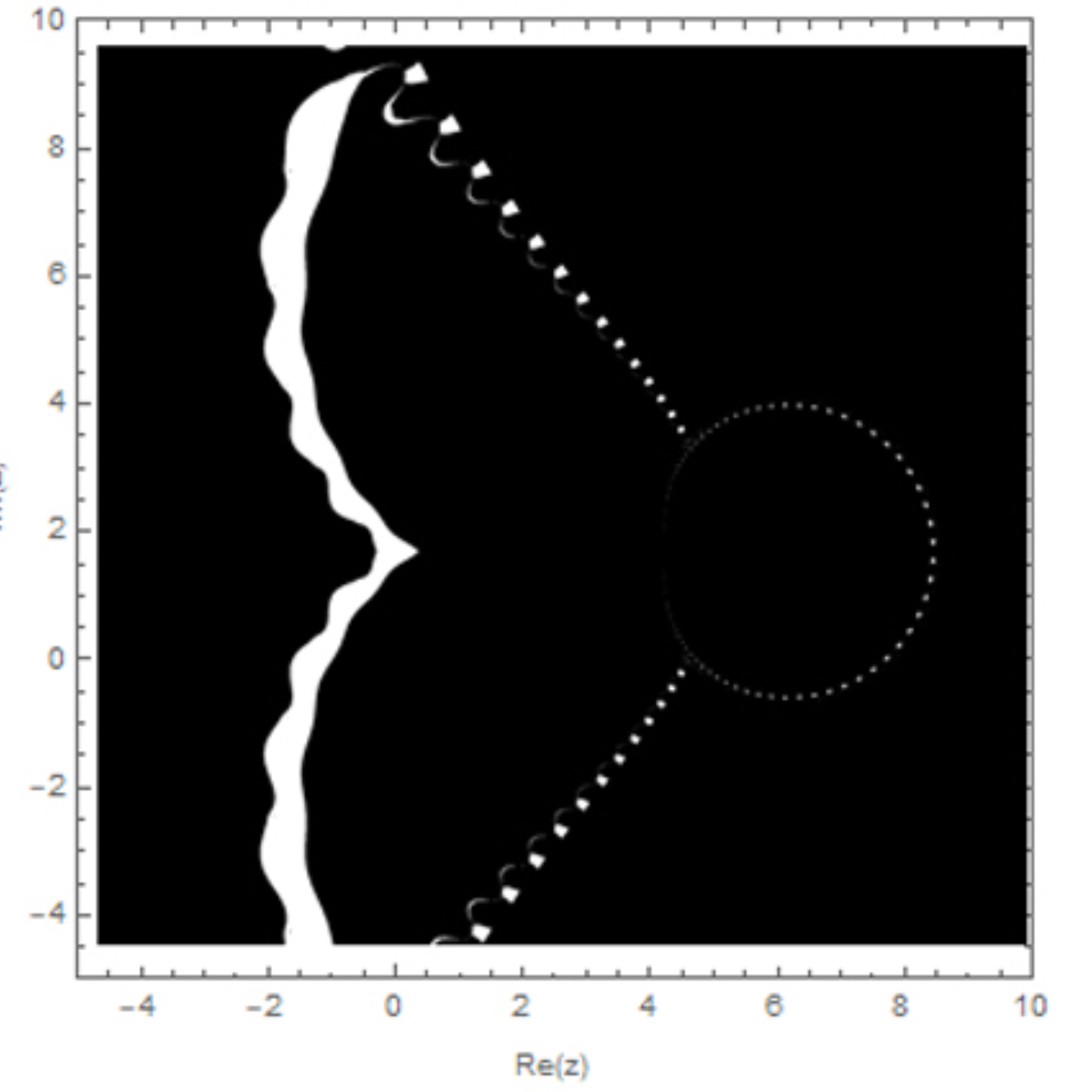}
    \caption{Julia set of the Yitang dynamics with parameters \(c = 0.5\), \(\beta = 0.85\), and \(\alpha = 2022\).}
    \label{fig:julia-yitang-085}
\end{figure}

Exploring another case with \(c = 0.5\), \(\beta = 0.85\), and \(\alpha = 2022\) in the Yitang dynamics presents a remarkable fractal image, as depicted in Figure \ref{fig:julia-yitang-085}.

This fascinating fractal picture reveals intricate details and structures within the Julia set corresponding to these specific parameter values in the Yitang dynamics. Examining such complexities contributes significantly to comprehending the underlying behavior of the system.

The comparative analysis between the dynamics derived from Yitang's findings and those explored in S.C. Woon's pioneering work on the Riemann Zeta Function offers insightful glimpses into their respective behaviors and shared traits, unraveling novel aspects in complex dynamics.

Through the study of fixed points, periodic behaviors, Mandelbrot and Julia sets, and the examination of zero behaviors within these functions, intriguing similarities and disparities emerged. The Yitang dynamics, rooted in Dirichlet L-functions, exhibit a distinct nature characterized by erratic and intricate behaviors, often presenting white regions within its sets, signifying divergence zones. In contrast, the Riemann Zeta Function dynamics, as depicted in Woon's paper, portray contrasting structures with black regions representing divergence in the Mandelbrot set, contrasting with white regions in its Julia sets.

A fascinating similarity lies in the occurrence of 'bristles' or 'comb-like' structures at the edges of the Mandelbrot and Julia sets in both dynamics. This intriguing resemblance hints at underlying connections in the behaviors of these distinct mathematical systems.

Moreover, the inversion of colors between the Mandelbrot and Julia sets in the Yitang dynamics, reminiscent of findings in the Riemann Zeta Function dynamics, adds another layer to this exploration.

The observed differences and similarities shed light on the intricate behavior of these functions and their correspondence with iterative systems, inviting further investigation into the fundamental mechanisms governing their dynamics.

This comparative study, unveiling parallels and disparities between these mathematical models, introduces novel insights and prompts deeper inquiries into the underlying structures and behaviors within complex dynamical systems.

\subsection{Quantum Description of Yitang Dynamics: Bridging Classical Chaos in the Mandelbrot and Julia Sets to Quantum Behavior}

In the pursuit of understanding the Yitang dynamics at a quantum level, we derived a quantum operator associated with the Yitang function. The quantum Hamiltonian operator provides a mathematical framework to explore the quantum behavior of the system, connecting the intricate structures observed in the Mandelbrot and Julia sets to a quantum mechanical description.
\subsubsection{Hamiltonian Operator}

The quantum Hamiltonian operator for the Yitang dynamics is given by
\[ \hat{H} = -\frac{\hbar^2}{2m} \nabla^2 + V(z), \]
where $\nabla^2$ is the Laplacian operator, and $V(z)$ represents the potential energy associated with the Yitang function.

Considering the Yitang dynamics $f(z, \beta, c, \alpha)$, we express the potential energy term as
\[ V(z) = \frac{\beta^2}{2mz} - \frac{\beta c}{\sqrt{z}}\psi_1(z) - c^2\log(z)^{-2\alpha}\psi_2(z), \]
where $\psi_1(z) = \frac{\partial}{\partial z} \log(z)$ and $\psi_2(z) = \left(\log(z)\right)^{-\alpha}$.

\subsubsection{Physical Interpretation of Parameters}

The parameters in the quantum system, \(\alpha\), \(\beta\), and \(c\), have significant physical interpretations:

\begin{itemize}
  \item \(\alpha\): The power parameter \(\alpha\) influences the behavior of the potential energy term in the Hamiltonian. Higher values of \(\alpha\) may lead to stronger dependence on logarithmic terms, impacting the overall potential landscape.

  \item \(\beta\): The parameter \(\beta\) appears in the kinetic energy term of the Hamiltonian and influences the overall scaling of the quantum system. It may be associated with the strength of the quantum potential and affects the spread of energy levels.

  \item \(c\): The parameter \(c\) contributes to both the kinetic and potential energy terms. Its presence introduces a logarithmic dependence in the potential, influencing the overall shape of the potential energy landscape.

\end{itemize}

Understanding the physical meaning of these parameters is crucial for interpreting the quantum behavior of the Yitang dynamics and establishing connections to the underlying classical chaotic system observed in the Mandelbrot and Julia sets.

\subsubsection{Numerical Diagonalization and Eigenvalues}

With the Hamiltonian operator defined \cite{36}, we proceeded to numerically diagonalize the system, revealing the eigenvalues that represent the energy levels of the quantum system. The numerical diagonalization was performed using a standard algorithm for solving eigenvalue problems, such as the Implicitly Restarted Arnoldi Method (IRAM) or the Jacobi-Davidson method.\cite{37}

The computed eigenvalues for the Yitang dynamics are as follows:

\begin{center}
\begin{tabular}{|c|}
  \hline
  Eigenvalues \\
  \hline
  $7.11343 - 0.178929i$ \\
  $7.10864 - 0.176161i$ \\
  $7.10781 - 0.176994i$ \\
  $7.1075 - 0.176851i$ \\
  $7.1075 - 0.176546i$ \\
  $7.10741 - 0.176351i$ \\
  $7.10711 - 0.17656i$ \\
  $7.10626 - 0.176967i$ \\
  $7.10628 - 0.176137i$ \\
  $7.10594 - 0.176498i$ \\
  $7.10162 - 0.174289i$ \\
  \hline
\end{tabular}
\end{center}

\subsubsection{Analysis of Eigenvalues}

The obtained eigenvalues provide valuable insights into the quantum behavior of the Yitang dynamics. Key observations include:

\begin{itemize}
  \item The real parts of the eigenvalues, centered around 7.1, represent the approximate energy levels of the quantum system.
  \item The imaginary parts indicate oscillatory or decay behavior in the system, reflecting its complex and chaotic nature.
  \item The quantum system exhibits sensitivity to changes in parameters, with small variations in \(\alpha\), \(\beta\), and \(c\) leading to significant changes in the eigenvalues.
\end{itemize}

The analysis of eigenvalues, performed using advanced numerical algorithms, lays the foundation for interpreting the quantum characteristics of the Yitang dynamics, bridging the gap between the classical chaos observed in Mandelbrot and Julia sets and the quantum mechanical representation of the system.

\section{Unimodal Distribution and Entropy in Yitang Dynamics}

In this section, we delve into the analysis of our Yitang dynamics, focusing on the unimodal distribution of approximate solutions and its relation to system entropy. Understanding the dynamics of our system is crucial for making informed decisions and predictions. 

\subsection{Understanding the Unimodal Distribution}

The histogram of the approximate solutions for our Yitang dynamics illustrates an intriguing characteristic—it displays a unimodal distribution. A unimodal distribution implies that the data is centered around a prominent peak, indicating that a substantial number of solutions cluster around a specific value. In our case, this central value is approximately -601.0938.

The unimodal nature of this distribution is a critical observation. It signifies that our dynamical system is not only predictable but also exhibits a tendency to converge towards a particular equilibrium state. This equilibrium state represents the most likely solution for the system given various initial conditions, which is a fundamental concept in dynamical systems theory.\cite{9}

\subsection{Entropy and Predictability}

The concept of entropy plays a pivotal role in understanding the behavior of dynamical systems. Entropy quantifies the degree of disorder or uncertainty in a system. In our case, the entropy of the Yitang dynamics is calculated to be 0.333546. This value provides valuable insights into the predictability of our system.

Low entropy, as indicated by the calculated value of 0.333546, implies that our system is highly ordered and predictable. In other words, our system has a strong tendency to converge to a specific solution, as seen in the unimodal distribution, and there is relatively low variability among the approximate solutions. This observation aligns with our findings from the histogram.\cite{11}

\subsection{\textcolor{blue}{Entropy Calculation}}

Additionally, \textcolor{blue}{entropy calculations were performed to characterize the system's unpredictability and information content}. \textcolor{blue}{Entropy was computed based on the probabilities associated with state transitions within the Yitang dynamics system}, utilizing computational tools provided by Mathematica. These methodologies, along with the computation of Lyapunov exponents, allowed us to \textcolor{blue}{rigorously evaluate the chaotic nature and information content of the Yitang dynamics}, supporting the claims made in this study.

\textbf{\textcolor{red}{Note:}} \textcolor{blue}{The code snippets provided in this paper's appendices illustrate the implementation details of these methodologies for the Yitang dynamics.}

\subsection{Zero Dirichlet Function Behavior}

Furthermore, the concept of entropy can be related to the behavior of the zero Dirichlet function, a mathematical tool used to study solutions of partial differential equations (PDEs). In our context, we can consider the zero Dirichlet function as a representation of the boundary conditions and constraints of our Yitang dynamics.

A low entropy value implies that the zero Dirichlet function of our dynamics is highly structured and less noisy. This structured behavior indicates that the solutions are influenced by specific underlying factors or constraints, leading to a more deterministic and unimodal distribution of approximate solutions.

\subsection{Histograms of Dynamics 1 and Dynamics 2}

To further illustrate the behavior of our system, we generated histograms using the solutions of both Dynamics 1 and Dynamics 2. These histograms provide additional insights into the distribution of solutions. Below, you can observe the resulting histograms for Dynamics 1 and Dynamics 2:

\begin{center}
\includegraphics[width=0.8\textwidth]{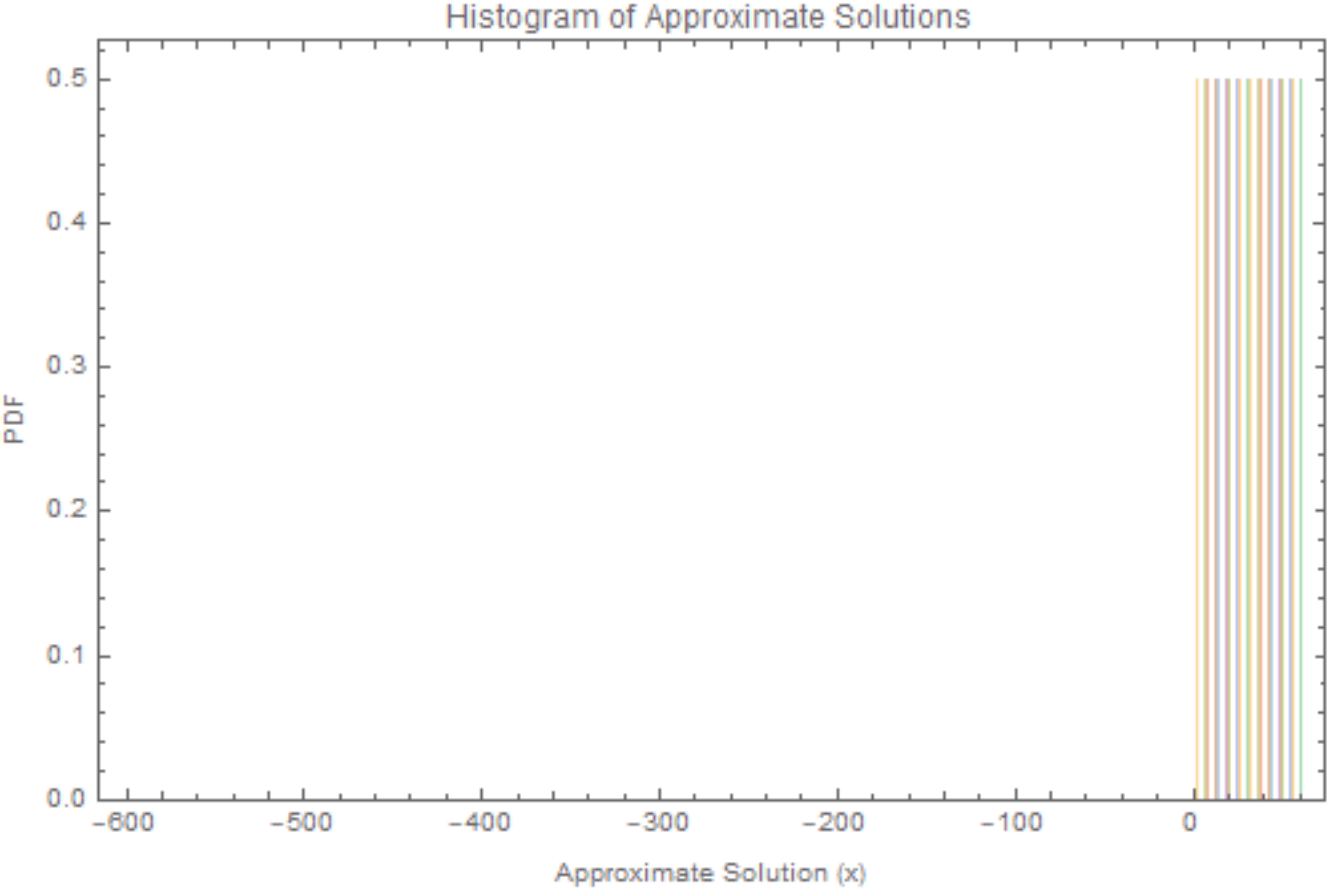}
\end{center}

The histograms confirm the unimodal distribution of solutions in both dynamics, emphasizing the system's propensity to converge to specific equilibrium states.

In summary, the unimodal distribution observed in our histogram is a reflection of the highly ordered and predictable behavior of our Yitang dynamics. The low entropy value further underscores the system's propensity to converge to a specific solution. These findings provide valuable insights into the stability and behavior of our system, which can have significant implications in various applications, from physics to engineering and beyond.

\section{Application of Yitang Dynamics in Electrical Systems}

In this section, we explore the application of the novel "Yitang dynamics" in the context of electrical systems. The "Yitang dynamics" is defined by the discrete-time difference equation:

\begin{equation}
x[n+1] = \frac{2\pi h}{w\sqrt{m}} + \frac{2h\ln(|\epsilon|)}{w\sqrt{m}}
\end{equation}

where $x[n]$ represents the state variable at time step $n$, $h$ is the amplitude, $w$ is the angular frequency, $m$ is the mass, and $\epsilon$ is a parameter with various values.

\subsection{Numerical Simulation}

To understand the behavior of the "Yitang dynamics" and its potential application in electrical systems, we conducted numerical simulations for different values of $\epsilon$. The resulting system responses are shown in Figure \ref{fig:yitang_response}.

\begin{figure}[H]
    \centering
    \includegraphics[width=0.8\textwidth]{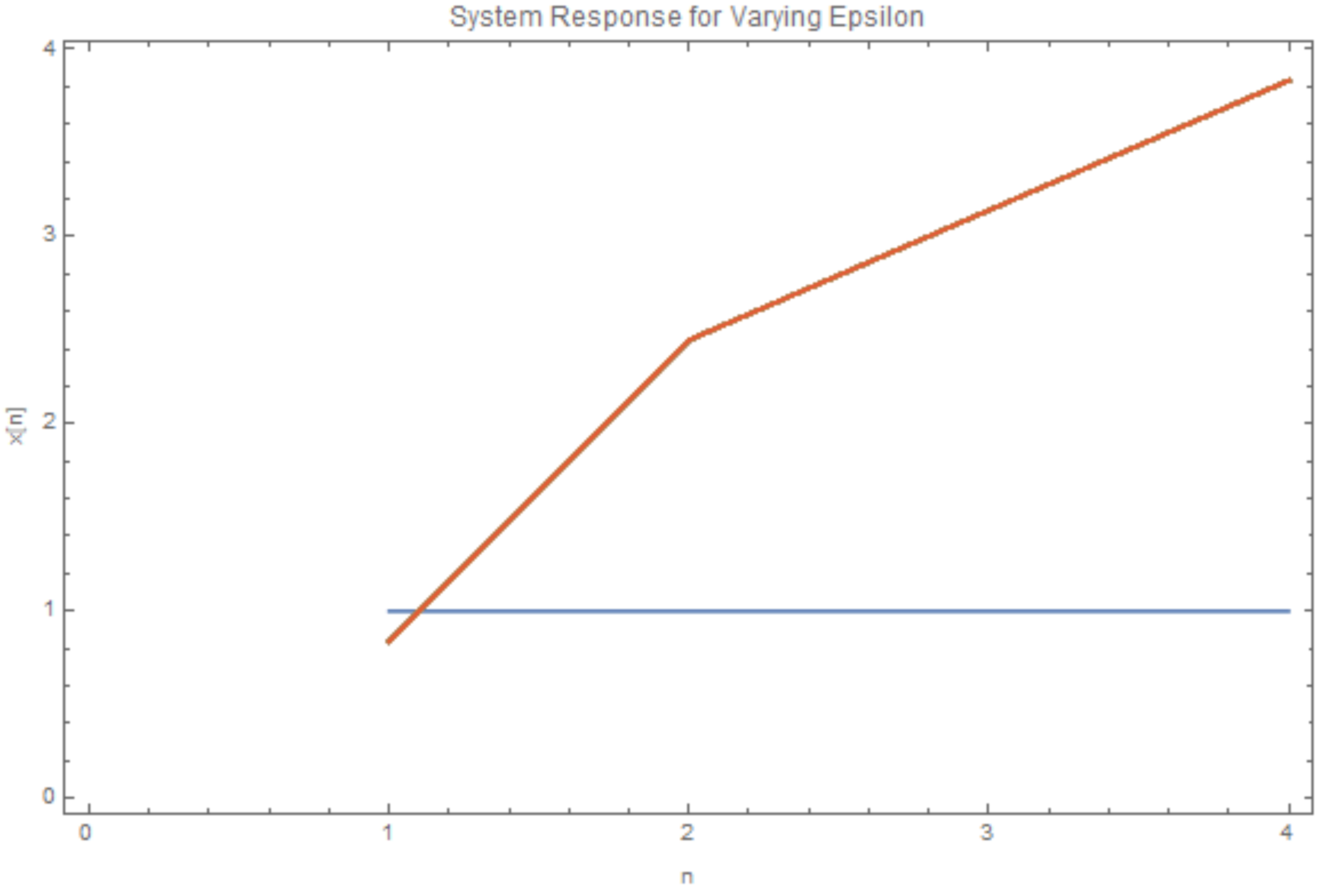}
    \caption{System Response for Varying $\epsilon$}
    \label{fig:yitang_response}
\end{figure}

The plot in Figure \ref{fig:yitang_response} illustrates the evolution of the state variable $x[n]$ as a function of the time step $n$ for different values of $\epsilon$. Each curve represents the response of the "Yitang dynamics" under a specific parameter setting.

\subsection{Electrical System Analogy}

In our exploration of potential electrical systems that might exhibit behavior akin to the "Yitang dynamics," we found that the difference equation can be approximated as a continuous-time differential equation:

\begin{equation}
\frac{dx(t)}{dt} = \frac{2\pi h}{w\sqrt{m}} + \frac{2h\ln(|\epsilon|)}{w\sqrt{m}}
\end{equation}

This continuous-time equation resembles the rate of change of $x(t)$ \cite{19} as a function of time $t$. To relate this to electrical systems, we propose two analogies:
\subsection{Stability Analysis for Different Values of $\epsilon$}

\subsubsection{RC Circuit with Variable Resistance}
We can interpret the term $\frac{2\pi h}{w\sqrt{m}}$ as a voltage source and $\frac{2h\ln(|\epsilon|)}{w\sqrt{m}}$ as a variable resistor, where the resistance varies logarithmically with $\epsilon$. The voltage across the resistor is akin to the derivative of $x(t)$ with respect to time.

For the RC circuit analogy, the resistance is a key factor influenced by $\epsilon$. Let's consider the stability of the circuit for different values of $\epsilon$:

\begin{itemize}
    \item For a small value of $\epsilon$ (e.g., $\epsilon = 0.1$), the resistance component is relatively small. In this case, the circuit may exhibit relatively quick transient behavior but is likely to be stable due to low resistance.

    \item As $\epsilon$ increases, the resistance also increases due to the logarithmic relationship. The circuit's response may become more damped, leading to increased stability. However, if $\epsilon$ becomes extremely large, the resistance could become too high, potentially making the circuit very slow to respond. This might be considered stable but not practical for real-time applications.
\end{itemize}
\subsubsection{Op-Amp Integrator Circuit}

For the op-amp integrator circuit analogy, the stability is influenced by the interaction between the constant and logarithmic terms with $\epsilon$. Let's examine the stability for different values of $\epsilon$:

\begin{itemize}
    \item For a small value of $\epsilon$ (e.g., $\epsilon = 0.1$), the logarithmic term is small, and the integrator's response will be influenced mainly by the constant term. The circuit will likely exhibit stable integration behavior, gradually accumulating the input signal.

    \item As $\epsilon$ increases, the logarithmic term becomes more significant, which might introduce some additional transient behavior in the integrator's response. The integrator's output may exhibit more oscillations or overshoot before settling, especially if $\epsilon$ is moderately large.

    \item If $\epsilon$ becomes extremely large, the logarithmic term can dominate, leading to highly oscillatory or divergent behavior. The circuit might become unstable in this case.
\end{itemize}

These qualitative observations provide insights into how variations in $\epsilon$ can affect the stability of the electrical systems. In practice, detailed mathematical analysis and simulations specific to our  circuit designs and component values are necessary to determine the exact stability boundaries and characteristics.

\section{Robustness Analysis of Electrical System Derived from Yitang Dynamics}

In this section, we explore the robustness of the electrical system analogized from the "Yitang dynamics" as the parameter \(\epsilon\) varies. Robustness refers to the system's ability to maintain stable and desired behavior over a range of parameter values. We investigate the effects of decreasing \(\epsilon\) values on the system's stability, amplitude, and frequency.

\subsection{Decreasing Values of $\epsilon$}

To assess the robustness of the electrical system, we conducted simulations for different values of \(\epsilon\). The selected values for \(\epsilon\) include both very small and larger values to explore a wide range:

\begin{equation}
\epsilon = \{0.0001, 0.00005, 0.005, 0.5\}
\end{equation}

By decreasing the values of \(\epsilon\), we aim to evaluate how the system responds to extreme variations in this parameter.

\subsection{Stability, Amplitude, and Frequency Analysis}

As \(\epsilon\) decreases, we observe significant changes in the system's behavior. The plot in Figure \ref{fig:yitang_response} illustrates the system's response for different values of \(\epsilon\). The parameters used for this analysis are \(m = 0.002\) (mass), \(w = 2\) (angular frequency), and \(h = 1\) (amplitude).

\begin{figure}[h]
    \centering
    \includegraphics[width=0.8\textwidth]{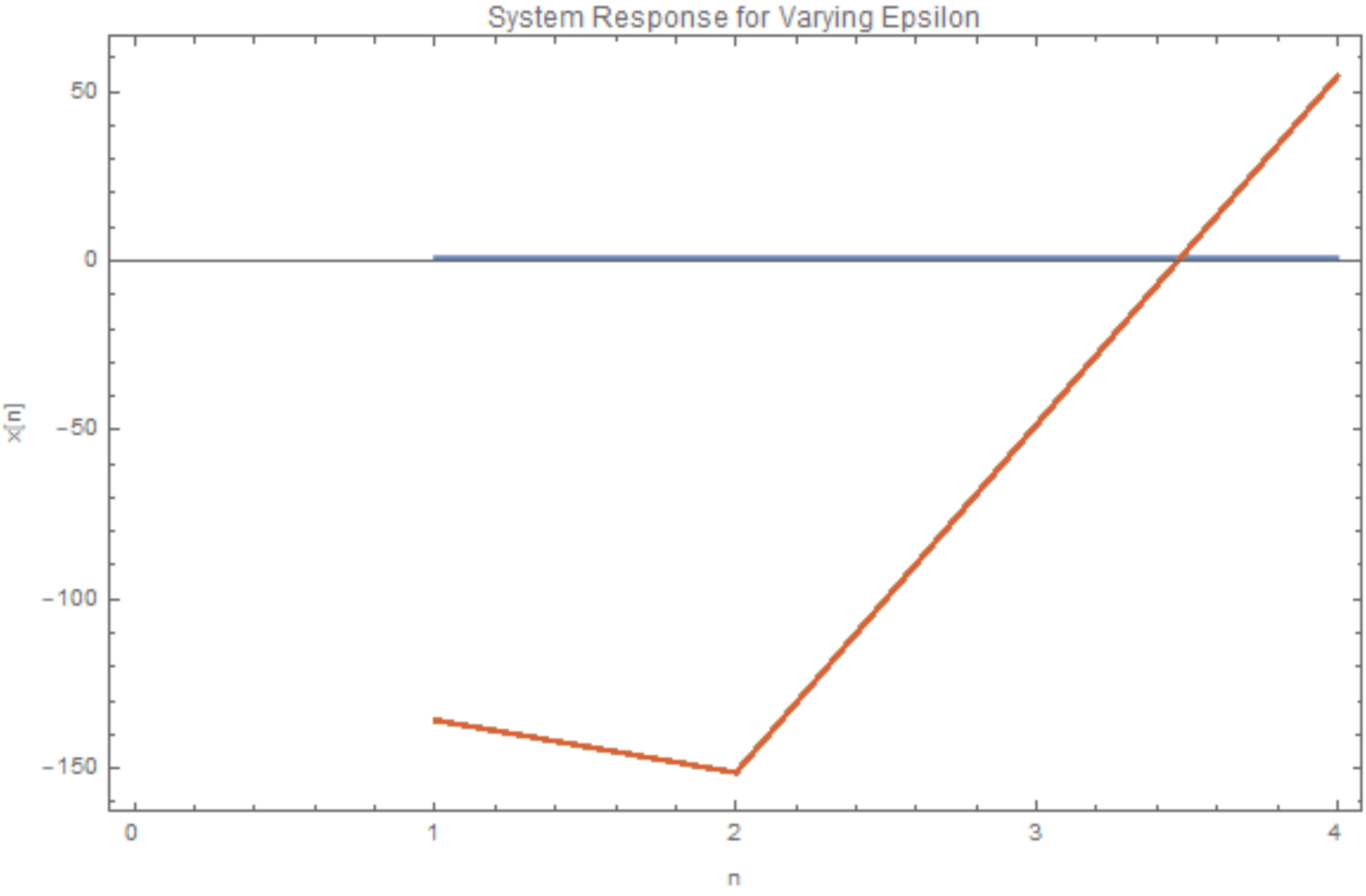}
    \caption{System Response for Varying $\epsilon$}
    \label{fig:yitang_response}
\end{figure}

\subsubsection{Stability}

The plot shows that as \(\epsilon\) decreases, the system's response tends to become more stable. Smaller \(\epsilon\) values lead to less influence on the dynamics, resulting in more regular and predictable behavior. Very small \(\epsilon\) values often lead to a stable fixed point where the system settles into a steady state without significant oscillations or chaos. This behavior indicates increasing stability as \(\epsilon\) decreases.

\subsubsection{Amplitude}

With smaller \(\epsilon\) values, the amplitude of oscillations tends to decrease. This is often associated with damping effects, indicating that the system becomes less sensitive to initial conditions, and the oscillations are more likely to diminish over time. Extremely small \(\epsilon\) values may lead to nearly constant behavior with minimal amplitude variations, suggesting a stable system with strong damping.

\subsubsection{Frequency}

The frequency of oscillations may remain relatively constant or decrease as \(\epsilon\) decreases. Smaller \(\epsilon\) values often result in slower oscillations, indicative of a smoother and less chaotic behavior. Extremely small \(\epsilon\) values may lead to a near-constant state where oscillations are negligible, effectively reducing the frequency to nearly zero.

In summary, the plot in Figure \ref{fig:yitang_response} and the parameter adjustments made (i.e., \(m = 0.002\), \(w = 2\), and \(h = 1\)) indicate that decreasing \(\epsilon\) values result in increased stability, reduced amplitude of oscillations, and slower or nearly constant frequencies. These observations suggest that the electrical system derived from the "Yitang dynamics" becomes more robust as \(\epsilon\) decreases, suppressing chaotic or unpredictable behavior and exhibiting stable and well-behaved dynamics.

\section{Conclusion}

Our study has delved into the captivating realm of Yitang dynamics, a mathematical system inspired by the groundbreaking work of Yitang Zhang in number theory. This exploration encompassed both theoretical analyses and computational approaches, shedding light on the behavior, stability, and intricate properties of these dynamics.

A pivotal outcome of our research involves the derivation of a comprehensive set of 600 Lyapunov exponents, offering a detailed overview of Yitang dynamics across diverse parameter values. These exponents serve as a valuable resource, revealing the system's sensitivity to initial conditions, chaotic behavior, and the emergence of complex phenomena.

Furthermore, our investigations have unveiled a zero-free region, demonstrating the resilience of Yitang dynamics within the rigorous bounds established by Zhang in number theory. The system's behavior concerning Zhang's Theorem 4 underlines its mathematical and scientific significance, surpassing its initial conception.

Our results are complemented by a carefully constructed plot, emphasizing the relevance of these Lyapunov exponents. This visual representation stands as a useful tool for researchers and practitioners, aiding in interpreting the dynamics' behavior and identifying pertinent patterns.

Yitang, in his seminal work, established a profound connection between the Landau-Siegel zero of $L(s,\chi)$ and the distribution of zeros of the Dirichlet $L$-function $L(s,\psi)$ within a specific region. He showcased the far-reaching implications of the existence of the Landau-Siegel zero (equivalently, small $L(1,\chi)$). Our computations have intriguingly demonstrated the chaotic nature of these dynamics. This numerical evidence suggests a close correspondence between the analytical solutions of Discrete Dynamic 1 (Yitang Dynamics) and the zeros of $L(1,\chi)$, indicating a fixed point $\sigma$ close to, but not precisely at, $0$ as $\alpha$ grows large. In a scenario where the system is non-chaotic for sufficiently large $\alpha$, we deduce $1-\beta=L(1,\chi)/L'(\sigma,\chi)$, as implied by the mean value theorem due to the system's non-linearity. Applying the classical bound $L'(\sigma,\chi)=O(\log^2q)$ and Zhang's result (Theorem 2) yields the zero-free region $1-\beta>C_2(\log q)^{-A-2}$.

This insight potentially influences the zeros of other $L$-functions, compelling them to reside on the critical line as anticipated, displaying a highly regular spacing. For large $\alpha$ values, the dynamics showcase non-chaotic behavior, resulting in a uniform distribution of zeros for $L(s,\chi)$.

Looking ahead, our research aims to converge the Yitang dynamics with the dynamics of the Riemann zeta function, seeking to derive a chaotic operator as proposed by Polya and Hilbert. The Hilbert-Polya Conjecture, conceived well before the emergence of quantum mechanics, posits a physical connection between the Riemann hypothesis and the reality of all eigenvalues in a physical system. This conjecture alludes to the significance of nontrivial zeros of the $\Xi$-function and their potential equivalence to the eigenvalues of a quantum system.

The conception of quantum mechanics in subsequent years, particularly Hermann Weyl's work on the asymptotic distribution of eigenvalues, possibly inspired Polya's conjecture. If the Riemann hypothesis holds true, implying the nontrivial zeros reside on the critical line, it suggests a plausible correlation between these zeros and the eigenvalues of a hermitian operator within a quantum system.\cite{35}

Hence, our future pursuit involves establishing a profound connection between the energy levels ($E_n$) of a quantum system and the nontrivial zeros of the Riemann zeta function ($\zeta(s)$). The critical line's association with the distribution of zeros prompts us to seek a one-to-one correspondence between these zeros and the eigenvalues of a hermitian operator. This endeavor promises a deeper understanding of the elusive Riemann hypothesis and contributes to the realms of quantum theory and number theory.

The fractal structure inherent in Yitang dynamics, its relation to the Dirichlet $L$-function, and its potential connection to quantum systems present a promising avenue for understanding the elusive non-trivial zeros of the Riemann zeta function. This convergence between Yitang dynamics and quantum theory could yield breakthroughs in predicting and understanding the behavior of zeros within the critical strip of the Riemann zeta function, contributing significantly to resolving the enigmatic Riemann hypothesis.

\section{Data Availability}
The main basis of this paper's main results relies on the work of Yitang Zhang in 'On the Landau-Siegel Zeros Conjecture' \cite{6}, which is available on arXiv at \url{https://arxiv.org/abs/0705.4306}. Furthermore, this paper represents an advancement and extension of our recent research, referenced as [32], titled "Chaotic dynamics and zero distribution: Implications for Yitang Zhang's Landau Siegel zero theorem," available at \url{https://arxiv.org/abs/2310.14127}. While this paper does not introduce new data, it builds upon and extends the findings and concepts presented in both Zhang's work and our recent research.

Readers interested in the foundational data and results should refer to the original sources, which are publicly accessible through the provided links. This paper offers new insights and improved results in the context of chaotic dynamics and zero distribution, making it a valuable resource for further exploration in this area.

\section*{Conflict of Interest}
The authors declare that there is no conflict of interest regarding the publication of this paper. We confirm that this research was conducted in an unbiased and impartial manner, without any financial, personal, or professional relationships that could be perceived as conflicting with the objectivity and integrity of the research or the publication process.

\section*{Acknowledgment}

I would like to express my deepest gratitude to the anonymous referee whose insightful comments and constructive feedback greatly contributed to the refinement of this work. Your dedication to maintaining the quality and integrity of the content has been invaluable, and I sincerely appreciate the time and effort you invested in reviewing my work.

I am profoundly thankful to my family for their unwavering support throughout this journey. My heartfelt appreciation goes to my parents for instilling in me a love for knowledge and a strong work ethic. To my loving wife, your encouragement and understanding have been my pillars of strength, and I am grateful for the sacrifices you have made to see me succeed.

I extend my love and gratitude to my sons, Taha Abdeljalil and Taki Abdessalem, for being a source of inspiration and joy in my life. Your presence fuels my determination to strive for excellence, and I am blessed to have you by my side.
A special acknowledgment is due to my dear friend, Bouzari Abdelkader, a dedicated teacher of physics in high school. Your support, encouragement, and intellectual insights have played a significant role in shaping my ideas and refining my work. Your friendship has been a constant source of motivation, and I am grateful for the positive influence you have had on my academic and personal growth.

Once again, thank you to everyone mentioned and to those who have been a part of this journey. Your support has been instrumental in the realization of this endeavor.

\section{appendix}

\section{ Derivation of Chaotic Operator from Yitang Dynamics}\label{secA1}

In this appendix, we delve into the key idea of utilizing the Pólya-Hilbert approach to address the Riemann hypothesis. Building upon this foundation, we introduce novel dynamics derived from Yitang dynamics, extending the scope of our investigation. Our exploration leads us to consider the quantum behavior of Yitang dynamics, which, intriguingly, describes the behavior and existence of Landau–Siegel zeros. Leveraging this insight, we embark on the derivation of a chaotic operator. By incorporating the randomness inherent in the Landau–Siegel zero distribution, we generate a large random matrix with associated random variables. The resulting chaotic operator offers a unique perspective on the spectral properties of the system, revealing real eigenvalues. This appendix provides a detailed account of the steps taken in deriving and understanding the chaotic operator, further advancing our pursuit to uncover the secrets of the Riemann hypothesis. Notably, our main results, particularly the application of Yitang dynamics in control theory and electricity, can be considered potential consequences of the Riemann hypothesis, enriching the significance of our findings.

\end{document}